\documentclass{m2an}
\usepackage{graphicx}

\usepackage{caption}
\usepackage[caption=false]{subfig}

\usepackage{mathtools}
\usepackage{amssymb}
\usepackage{dsfont}
\usepackage{booktabs}
\usepackage{bm}
\usepackage{color}
\usepackage[mathcal]{eucal}
\usepackage{pgfplots}
	\pgfplotsset{compat=1.3}
\usepgfplotslibrary{groupplots}
\usetikzlibrary{matrix}

\newtheorem{Assumption}[thrm]{\sc Assumption}
\def\norm#1{\left|\!\left| #1 \right|\!\right|}
\def\fp{\,{:}\,}

\def\tens#1{\pmb{\mathsf{#1}}}
\def\vec#1{\boldsymbol{#1}}
\def\R{\mathbb{R}}
\def\Rdd{{\mathbb{R}}^{d\times d}}
\def\sym{\mathop{\mathrm{sym}}\nolimits}
\def\tr{\mathop{\mathrm{tr}}\nolimits}
\def\co{\mathop{\mathrm{co}}\nolimits}
\def\Rds{\mathbb{R}^{d \times d}_{\sym}}

\def\diver{\mathop{\mathrm{div}}\nolimits} 
\def\Du{\BD(\bu)}
\def\wconv{\rightharpoonup}
\newcommand{\Lebs}[4][\Omega]{L^{#2}_{#3}(#1)^{#4}}
\newcommand{\Lp}[2][\Omega]{\Lebs[#1]{#2}{}{}}

\newcommand{\Lmean}[2][\Omega]{\Lebs[#1]{#2}{0}{}}
\newcommand{\Lsym}[2][\Omega]{\Lebs[#1]{#2}{\sym}{d \times d}}

\newcommand{\Lsymtr}[2][\Omega]{\Lebs[#1]{#2}{\sym,\tr}{d \times d}}
\newcommand{\Sobs}[5][\Omega]{W^{#2, #3}_{#4}(#1)^{#5}}

\def\bh{\vec{h}}
\def\bu{\vec{u}}
\def\bv{\vec{v}}
\def\BD{\tens{D}}
\def\BG{\tens{G}}
\def\BH{\tens{H}}
\def\BI{\tens{I}}
\def\BS{\tens{S}}
\def\BT{\tens{T}}
\def\Bphi{\tens{\phi}}
\def\Bsigma{\tens{\sigma}}
\def\Btau{\tens{\tau}}
\def\Scr{\mathcal{{S}}}
\def\Vndiv{V^n_{\diver}}
\def\Signsym{\Sigma^n_{\sym}}
\def\Signtr{\Sigma^n_{\tr}}

\newcommand{\txtbl}[1]{\textcolor{black}{#1}}
\begin{document}
\title{A semismooth Newton method for implicitly constituted non-Newtonian fluids and its application to the numerical approximation of Bingham flow}\thanks{This work was supported by the Alexander von Humboldt-Stiftung; Chair in Applied Analysis (Alexander von Humboldt-Professur)}
\author{Pablo Alexei Gazca-Orozco}\address{Department of Mathematics, FAU Erlangen-N\"{u}rnberg, 91058 Erlangen, Germany; e-mail: \texttt{alexei.gazca@math.fau.de}
}
%
\date{\today}
\begin{abstract}  We propose a semismooth Newton method for non-Newtonian models of incompressible flow where the constitutive relation between the shear stress and the symmetric velocity gradient is given implicitly; this class of constitutive relations captures for instance the models of Bingham and Herschel--Bulkley. The proposed method avoids the use of variational inequalities and is based on a particularly simple regularisation for which the (weak) convergence of the approximate stresses is known to hold. The system is analysed at the function space level and results in mesh-independent behaviour of the nonlinear iterations.
\end{abstract}
%
%
\maketitle
\section{Introduction}
\label{sec;introduction}

Let $\Omega\subset \mathbb{R}^d$ be a Lipschitz polygonal/polyhedral domain for $d\in \{2,3\}$ and consider the following system for the velocity $\bu \colon \Omega \to \R^d$, shear stress tensor $\BS \colon \Omega \to \Rds$ and pressure $p \colon \Omega \to \R$ of an incompressible fluid:
\begin{subequations}\label{eq:PDE}
\begin{align}
	\alpha \bu - \diver\BS + \diver(\bu&\otimes\bu) + \nabla p = \bm{f} && \text{in }\Omega,\notag\\
	\diver \bu &= 0 && \text{in }\Omega,\label{eq:PDE_}\\
\bu &= \bm{0} && \text{ on }\partial\Omega,\notag
\end{align}
where $\bm{f}\colon \Omega \to \R^d$ is given, and $\alpha\geq 0$ is e.g.\ a parameter that arises from an implicit time discretisation ($\alpha = 0$ for the steady problem). In this work we will close the system above with an implicit constitutive relation of the form
\begin{equation}
	\BG(\BS,\Du) = \bm{0} \qquad\qquad \text{a.e. in }\Omega,\label{eq:PDE_CR}\\
\end{equation}
\end{subequations}
where $\Du := \frac{1}{2}(\nabla\bu + \nabla \bu^\top)$ is the symmetric velocity gradient and $\BG\colon \Rds\times \Rds \to \Rds$ defines a monotone graph on the space of symmetric matrices (the precise assumptions on the function $\BG$ will be introduced later). The theory of implicitly constituted fluids, introduced in \cite{Rajagopal2006,Rajagopal:2003}, provides a framework that captures a wide range of models and allows for a thermodynamically consistent analysis (see e.g.\ \cite{Rajagopal2008}). For example, the usual Newtonian and power-law models can be written in the form \eqref{eq:PDE_CR}:
\begin{subequations}\label{eq:Newtonian_Powerlaw}
\begin{align}
	\BG(\BS,\Du) &:= \BS - 2\nu\Du &&\qquad \nu >0,\\
	\BG(\BS,\Du) &:= \BS - K|\Du|^{r-2}\Du && \qquad K>0,\, r>1.
\end{align}
\end{subequations}

As a motivating example, in this work we will place special emphasis on the Bingham constitutive relation for viscoplastic fluids, defined by the dichotomy:
\begin{equation}\label{eq:Bingham_intro}
  \renewcommand{\arraystretch}{1.2}
\left\{
	\begin{array}{ccc}
		|\BS|\leq \tau_* & \Longleftrightarrow & \Du = \bm{0}, \\
|\BS|> \tau_* & \Longleftrightarrow & \BS = 2\nu\Du + \tau_*\displaystyle\frac{\Du}{|\Du|},
	\end{array}
\right.
\end{equation}
where $\tau_*\geq 0$ is the yield stress, and $\nu>0$ is the viscosity; the closely related Herschel--Bulkley constitutive relation can be obtained e.g.\ by letting $\nu = \nu(\Du):= K|\Du|^{r-2}$, with $K>0$ and $r>1$. Such models can be used to describe waxy crude oil, toothpaste, paint, pastes, drilling muds, and mango jam, among other things \cite{Bird1983,Glowinski2010,Basu,Luu2009}; see \cite{Balmforth2014} for a nice survey on various aspects of the modeling and simulation of viscoplastic fluids. Note that, in contrast to the constitutive relations \eqref{eq:Newtonian_Powerlaw}, the  relation \eqref{eq:Bingham_intro} cannot be written  in terms of a (single-valued) function $\BS = \mathcal{S}(\Du)$. However, within the framework of implicitly constituted fluids it can be very naturally written, for instance, using either of the two following expressions:
\begin{subequations}\label{eq:Bingham_implicit}
\begin{gather}
	\BG(\BS,\Du) := (|\BS| - \tau_*)^+\BS
	-  2\nu(\tau_* + (|\BS|-\tau_*)^+)\Du = \bm{0}, \label{eq:Bingham_implicit_1} \\
	\BG(\BS,\Du) := |\Du|\BS - (\tau_* + 2\nu |\Du|)\Du = \bm{0},\label{eq:Bingham_implicit_2}
\end{gather}
\end{subequations}
where $p^+ := \max\{p,0\}$ is the positive part of $p$. For a rigorous mathematical analysis of models describing implicitly constituted incompressible fluids, the reader is referred to \cite{Bulicek:2009,Bulicek:2012}; in particular existence of weak solutions for large data is known to hold both for the steady and unsteady problems, under appropriate assumptions for $\BG$ (see also \cite{Blechta2019} for a classification of various models that can be described within this framework). Regarding the numerical analysis, in \cite{Diening:2013} and \cite{Sueli2018} the authors established (weak) convergence of the finite element approximations to a weak solution of the system for the steady and unsteady problems, respectively. This was later extended to formulations including the shear stress $\BS$ in \cite{Farrell2019}, and an augmented Lagrangian preconditioner was proposed in \cite{Farrell2020a}.

In previous works dealing with the numerical analysis of implicitly constituted fluids a fixed point argument is employed to obtain existence of discrete solutions, and then the emphasis is placed on proving that convergence to a weak solution of the problem holds. However, there is barely any discussion on how to solve such discrete nonlinear systems. For smooth constitutive relations, such as \eqref{eq:Newtonian_Powerlaw}, this is not an issue since the standard Newton's method can be applied, but it is an important matter when dealing e.g.\ with the Bingham constitutive relation. In practice, it is common to circumvent the non-differentiability by regularising; a popular regularisation is for example the one proposed in \cite{Bercovier1980}:
\begin{equation}\label{eq:BE_regularisation}
	\BS_\varepsilon = \Scr_\varepsilon(\BD(\bu_\varepsilon))
	:= 2\nu \BD(\bu_\varepsilon) + \tau_*\frac{\BD(\bu_\varepsilon)}{\sqrt{\varepsilon^2 + |\BD(\bu_\varepsilon)|^2}}
	\qquad \varepsilon>0.
\end{equation}
Regularisations such as \eqref{eq:BE_regularisation} are very popular owing to their simplicity. However, a major drawback of this approach is that \txtbl{by and large} there are no results that guarantee the convergence of the shear stress $\BS_\varepsilon \to \BS$, as $\varepsilon\to 0$, and previous works suggest that this convergence might not actually hold in general (see e.g.\ \cite{Saramito2017,Frigaard2005,Putz2009}); \txtbl{one exception is the regularisation of the Bingham model introduced in \cite{DelosReyes2010}, for which the authors prove convergence of the approximate solutions to a solution of the original problem}. The main alternative, whenever an accurate approximation of the yield surfaces (i.e.\ the region of transition between the two regimes described by \eqref{eq:Bingham_intro}) is essential, is given by the augmented Lagrangian method, which is based on a variational inequality formulation of the problem (see \cite{Glowinski2010}). The main drawback of this method is that it requires the use of more sophisticated mathematical tools to be analysed and convergence can be slow; the search for ways to accelerate the algorithm is still a very active area of research (see e.g.\ \cite{Dimakopoulos2018}). 

Here we propose the use of an alternative regularisation, introduced in \cite{Bulicek2020}, for which the convergence $\BS_\varepsilon \wconv \BS$ can in fact be established. This regularisation takes the very simple form
{\color{black}
\begin{equation}\label{eq:CR_regularised}
\BG^{\varepsilon_1}_{\varepsilon_2}(\BS,\BD(\bu)) :=
	\BG(\BS - \varepsilon_1\BD(\bu),
	\BD(\bu) - \varepsilon_2\BS)\qquad
	\varepsilon_1,\varepsilon_2>0.
\end{equation}
}
The approach proposed here is conceptually very simple and avoids the use of formulations based on variational inequalities. Note that the regularised function $\BG_{\varepsilon_2}^{\varepsilon_1}$ inherits the smoothness of the original constitutive relation, and so in general it will not be differentiable (it brings however other advantages that will be discussed later). For this reason a semismooth version of Newton's method will be required.

In the classical version of Newton's method for solving the equation $F(z) = 0$, where $F\colon Z \to X$ is a function defined between two Banach spaces $Z$ and $X$, one applies iteratively the following steps until a convergence criterion is satisfied:
\begin{enumerate}
	\item Solve the linear problem: $DF(z^k)d^k = -F(z^k)$;
	\item Update the solution: $z^{k+1}:= z^k + d^k$.
\end{enumerate}
With the appropriate hypotheses it is well known that the iterates produced by this method converge (locally) quadratically to the solution of the problem. In the semismooth version of the algorithm one has to replace $DF(z^k)$ in step $1$ above by an element of the generalised Jacobian of $F$. Then, assuming some regularity conditions, it is possible to prove that the iterates converge at least (locally) superlinearly to the solution (the required hypotheses will be described in detail in subsequent sections).

After introducing the relevant weak and finite element formulations in Section \ref{sec:WeakandFE_formulations}, we will proceed in Section \ref{sec:semismoothNewton} to write the problem \eqref{eq:PDE_} with the regularised constitutive relation \eqref{eq:CR_regularised} in terms of a function $F$ between appropriate Banach spaces and formulate a semismooth Newton algorithm that, with appropriate assumptions, can be shown to fall within the framework of semismooth Newton methods in function spaces of \cite{Ulbrich2003}; in particular, we will see that some expressions describing the Bingham constitutive relation are more advantageous. The analysis is carried out at the function space level in order to avoid mesh-dependent behaviour. Numerical experiments focusing on the Bingham constitutive relation that support this result are carried out in Section \ref{sec:Examples}; these experiments employ the finite element method and a stress-velocity-pressure formulation.

Semismooth Newton methods have been applied before in the context of viscoplastic flow; semismooth Newton methods were analysed in \cite{DelosReyes2010,DelosReyes2011} for the solution of certain variational inequalities related to the description of Bingham flow, and so the approach is distinct to the one presented here (as mentioned above, one of the advantages the framework presented here is precisely that it avoids the formalism of variational inequalities). A non-smooth Newton method was developed for a different formulation describing Herschel--Bulkley fluids with $r\in (1,2)$ in \cite{Saramito2016}, but, as opposed to the present work, some mesh-dependence was observed in the solver. In addition, the analysis from \cite{Saramito2016} is closely tied to the particular structure of the Herschel--Bulkley relation, which is a constraint not present here: the framework employed here captures a wider variety of models; for instance, given the implicit nature of the constitutive relation, one could very easily consider stress-dependent viscosities or even swap the roles of the shear stress $\BS$ and the symmetric velocity gradient $\Du$, and obtain models describing inviscid (i.e.\ Euler) fluids for low values of $|\Du|$, and viscous otherwise (see \cite{Blechta2019} for a more detailed description of related models).

\section{Weak and finite element formulations}\label{sec:WeakandFE_formulations}
\subsection{Implicitly constituted fluids}
Throughout this work we will employ standard notation for Sobolev and Lebesgue spaces (for example $(W^{1,r}(\Omega), \norm{\cdot}_{W^{1,r}(\Omega)})$). We define the space $W^{k,r}_0(\Omega)$, with $k\in\mathbb{N}$, $r>1$, to be the closure of the space of compactly supported smooth functions $C^\infty_0(\Omega)$ with respect to the norm $\norm{\cdot}_{W^{k,r}(\Omega)}$. In addition, for $r\in (1,\infty)$ we define the following subspaces:
\begin{gather*}
L^r_0(\Omega) := \left\{q\in L^r(\Omega) \, : \, \int_\Omega q = 0\right\},\\
W^{1,r}_{0,{{\diver}}}(\Omega)^d := \overline{\{\bv\in C^\infty_0(\Omega)^d\, :\, {{\diver}}\,\bv=0\}}^{\|\cdot\|_{W^{1,r}(\Omega)}},\\
L_{{{\tr}}}^r(Q)^{d\times d} := \{\bm{\tau}\in L^r(Q)^{d\times d}\, : \, {{\tr}}(\bm{\tau}) = 0\},\\
L_{{\sym}}^r(Q)^{d\times d} := \{\bm{\tau}\in L^r(Q)^{d\times d}\, : \, \bm{\tau}^\top = \bm{\tau}\},\\
L_{{\sym,\tr}}^r(Q)^{d\times d} := L_{{\sym}}^r(Q)^{d\times d} \cap L_{{{\tr}}}^r(Q)^{d\times d}.
\end{gather*}
The operator $\tr(\Btau)$ above denotes the usual trace of the matrix function $\Btau$. The Frobenius inner product between matrices will be written as $\Btau\fp\Bsigma := \tr(\Btau\Bsigma^\top)$.

Suppose that $\bm{f}\in L^{r'}(\Omega)^d$. In order to focus solely on the difficulties brought by the nonlinear constitutive relation, we will neglect the convective term in the subsequent analysis. Therefore, in the weak formulation of the problem we look for $(\BS,\bu,p)\in  \Lsymtr{r'}\times V^r \times L^{r'}_0(\Omega)$, where $V^r := W^{1,r}_{0,\diver}(\Omega)^d \cap L^2(\Omega)^d$ and $r'$ is the H\"{o}lder conjugate of $r$,  such that

\begin{subequations}\label{eq:WeakFormulation}
\begin{alignat}{2}
\alpha\int_\Omega \bu\cdot \bv +
\int_\Omega \BS\fp  \BD(\bv)  &-\int_\Omega p  \diver\bv = \int_\Omega \bm{f}\cdot\bv\quad & \forall\,\bv\in C^\infty_0(\Omega)^d,\\
\BG(\BS&, \Du) = \bm{0} & \textrm{a.e. in }\Omega,\\
-\int_\Omega &q  \diver\bu = 0 & \forall\, q\in C^\infty_0(\Omega).
\end{alignat}
\end{subequations}

In the existence analysis of systems such as \eqref{eq:WeakFormulation}, the assumptions on the constitutive relation are traditionally written in terms not of the function $\BG$, but of its induced graph on $\Rds\times \Rds$. In this work we will prefer to work directly with the function $\BG$, following the approach first proposed in \cite{Bulicek2020}, which leads to conditions that are easier to check in practice.

\begin{Assumption}[{Properties of $\BG$}]\label{As:G}
We have for some $r>1$:
\begin{enumerate}
	\item [(G1)] \label{itm:G-1} $\BG$ is locally Lipschitz and $\BG(\bm{0},\bm{0})=\bm{0}$;
	\item [(G2)] \label{itm:G-2} We have for any $\Btau,\Bsigma\in\Rds$:
\begin{gather*}
\partial_{\Bsigma}\BG(\Bsigma,\Btau)\geq \bm{0} \qquad \partial_{\Btau}\BG(\Bsigma,\Btau)\leq \bm{0},    \\
 \partial_{\Bsigma}\BG(\Bsigma,\Btau)-\partial_{\Btau}\BG(\Bsigma,\Btau)>0 \quad  \partial_{\Btau}(\Bsigma,\Btau)(\partial_{\Bsigma}\BG(\Bsigma,\Btau))^\top \leq 0;
\end{gather*}
\item [(G3)] \label{itm:G-3} Either
\begin{equation*}
    \forall\, \Btau\in\Rds\quad \liminf_{\Bsigma\to\infty}\BG(\Bsigma,\Btau)\fp\Bsigma >0,
\end{equation*}
or
\begin{equation*}
    \forall\, \Bsigma\in\Rds\quad \liminf_{\Btau\to\infty}\BG(\Bsigma,\Btau)\fp\Btau <0;
\end{equation*}
\item [(G4)] \label{itm:G-4} For any $(\Bsigma,\Btau)\in\Rds\times \Rds$ such that $\BG(\Bsigma,\Btau)=\bm{0}$:
\begin{equation*}
    \Bsigma\fp \Btau \geq c_1(|\Bsigma|^{r'} + |\Btau|^r) - c_2,
\end{equation*}
for some $c_1,c_2>0$.
\end{enumerate}
\end{Assumption}

The assumptions (G1)--(G4) imply that $\BG$ induces a maximal monotone $r$-coercive graph on $\Rds\times \Rds$ (see \cite[Def. 3.1]{Bulicek2020}) and so they guarantee the existence of a weak solution to the problem \cite[Thm. 2.2]{Bulicek2020}. It can be seen that the Bingham relation \eqref{eq:Bingham_implicit} and relations of power-law type like \eqref{eq:Newtonian_Powerlaw} satisfy the assumptions above (see e.g.\ \cite{Bulicek2020,Tscherpel2018}).

\begin{figure}
\centering
\subfloat[C][{\centering Herschel--Bulkley relation ($r=1.7$).}]{{%
	\includegraphics[width=0.5\textwidth]{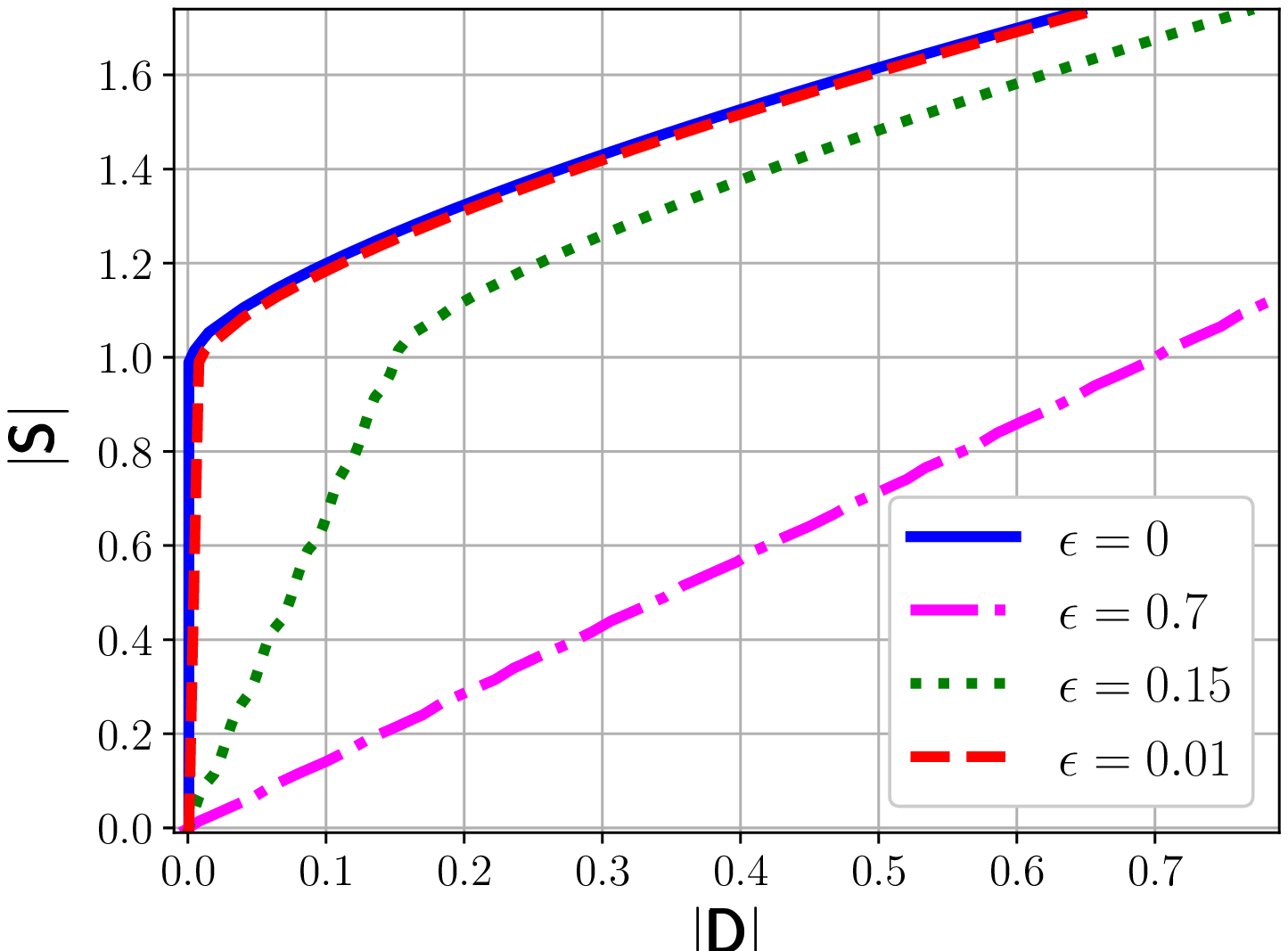}%
	}}%
	\subfloat[D][{\centering Bingham relation ($r=2.0$).}]{{%
	\includegraphics[width=0.5\textwidth]{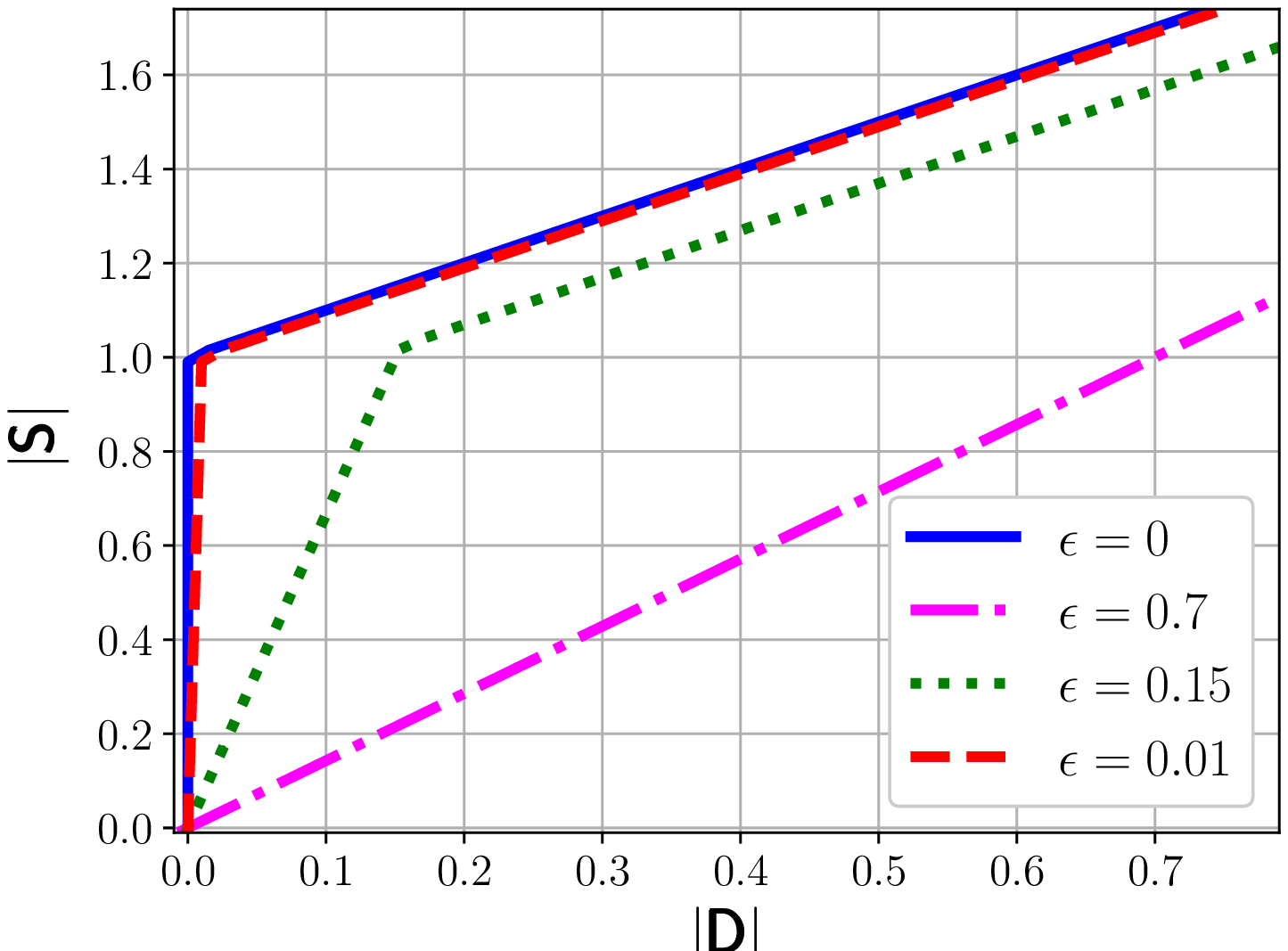}%
	}}%
	\caption{Regularisation \eqref{eq:CR_regularised} for the Bingham and Herschel--Bulkley constitutive relations with $\tau_*=1$, \textcolor{black}{taking the same value $\varepsilon=\varepsilon_1=\varepsilon_2$}.}%
	\label{fig:regularised_CR}
\end{figure}

The existence result obtained in \cite{Bulicek2020} relies on an approximation scheme based on solving the problem associated to the approximate constitutive relation \eqref{eq:CR_regularised} and then passing to the limit \textcolor{black}{$(\varepsilon_1,\varepsilon_2)\to 0$}; in particular, one has that the approximate solutions \textcolor{black}{$(\BS^{\varepsilon_1,\varepsilon_2},\bu^{\varepsilon_1,\varepsilon_2},p^{\varepsilon_1,\varepsilon_2})$} converge weakly to a weak solution $(\BS,\bu,p)$ of the system. From a theoretical point of view the regularisation \eqref{eq:CR_regularised} is advantageous because the resulting approximate graph is strictly monotone and $2$-coercive, thus leading to an approximate problem with a Hilbert space structure. Figure \ref{fig:regularised_CR} shows a plot of the regularised constitutive relation \eqref{eq:CR_regularised} for the Bingham and Herschel--Bulkley models for different values of $\varepsilon = \varepsilon_1 = \varepsilon_2$.

{\color{black}
\begin{rmrk}
	In the theoretical analysis carried out in \cite{Bulicek2020}, the parameters $\varepsilon_1$ and $\varepsilon_2$ were taken to be equal:
	\begin{equation}\label{eq:CR_regularised_sameeps}
		\BG_\varepsilon(\BS,\Du) :=
		\BG^{\varepsilon}_{\varepsilon}(\BS,\Du)=
		\BG(\BS - \varepsilon \Du, \Du - \varepsilon\BS) = 0 \qquad \varepsilon>0.
	\end{equation}
	This choice was made because having two distinct parameters did not bring any advantages to the analysis \cite[Rmrk.\ 4.2]{Bulicek2020}. Here we prefer to introduce the regularisation in the form \eqref{eq:CR_regularised} to emphasise the fact that the two parameters have different units: $\varepsilon_1$ has units of viscosity and $\varepsilon_2$ has units of fluidity; otherwise, starting directly from \eqref{eq:CR_regularised_sameeps} can give the impression that one is commiting a ``dimensional crime''. Furthermore, when performing numerical computations it is not clear whether different choices for $\varepsilon_1$ and $\varepsilon_2$ could result in better behaviour. However, for the sake of simplicity we will employ the form \eqref{eq:CR_regularised_sameeps} in the rest of this work.
\end{rmrk}
}

\begin{lmm}[{\cite{Bulicek2020}}]
	Suppose that $\BG$ satisfies Assumptions (G1)--(G4). Then the graph induced by $\BG_\varepsilon$ (recall \eqref{eq:CR_regularised_sameeps}) is uniformly monotone and $2$-coercive, i.e.\ for any $(\Bsigma_i,\Btau_i)\in \Rds\times \Rds$ such that $\BG_\varepsilon(\Bsigma_i,\Btau_i)= \bm{0}$, with $i\in\{1,2\}$, we have that
\begin{equation}\label{eq:coercivity_monotonicity_regularised}
\begin{gathered}
	\Bsigma_1\fp\Btau_1 \geq c(|\Bsigma_1|^2 + |\Btau_1|^2) - \tilde{c},\\
	(\Bsigma_1 - \Bsigma_2)\fp (\Btau_1 - \Btau_2)\geq c_\varepsilon (|\Bsigma_1- \Bsigma_2|^2 + |\Btau_1 - \Btau_2|^2),
\end{gathered}
\end{equation}
where $c,\tilde{c}>0$ and $c_\varepsilon = \frac{\varepsilon}{1 + \varepsilon^2}>0$.
\end{lmm}

\subsection{Finite element approximation}\label{sec:FEM}
We will now introduce the necessary notation to be able to define the finite element formulation of the system. Since the constitutive relation is implicit, a natural approach is to work with a 3-field formulation that includes the shear stress $\BS$ (such a formulation was analysed in \cite{Farrell2019}). Let $\{\mathcal{T}_n\}_{n\in\mathbb{N}}$ be a family of shape-regular triangulations such that the mesh size $h_n:= \max_{K\in \mathcal{T}_n}h_K$ vanishes as $n\to\infty$. We then define the following conforming families of finite element spaces:
\begin{align*}
	\Sigma^n &:= \left\{\Bsigma \in \Lsym{\infty} \colon \Bsigma|_K \in \mathbb{P}_{\mathbb{S}}(K)^{d\times d},\, K\in\mathcal{T}_n \right\}, \\
V^n &:= \left\{\bv \in \Sobs{1}{\infty}{0}{d} \colon \bv|_K \in \mathbb{P}_{\mathbb{V}}(K)^d,\, K\in\mathcal{T}_n \right\}, \\
M^n &:= \left\{q \in \Lp{\infty} \colon q|_K \in \mathbb{P}_{\mathbb{M}}(K),\, K\in\mathcal{T}_n \right\},
\end{align*}
where $\mathbb{P}_{\mathbb{S}}(K)$, $\mathbb{P}_{\mathbb{V}}(K)$, and $\mathbb{P}_{\mathbb{M}}(K)$ are spaces of polynomials defined on the element $K\in \mathcal{T}_n$. We will also introduce the following useful subspaces:
\begin{gather*}
	M_0^n := M^n \cap \Lmean{2},
	\quad \Signtr := \Sigma^n \cap L^2_{\sym,\mathrm{tr}}(\Omega)^{d\times d}, \\
	\Vndiv := \left\{\bv \in V^n \colon \int_\Omega q\diver \bv = 0\quad \forall q\in M^n \right\}.
\end{gather*}

\begin{Assumption}[Approximability]\label{as:approximability}
For every $s\in [1,\infty)$ we have that
\begin{align*}
	\inf_{\overline{\bm{v}}\in V^n}\|\bm{v}-\overline{\bm{v}}\|_{W^{1,s}(\Omega)} &\rightarrow 0 \quad \text{ as }n\rightarrow\infty\quad \forall\,\bm{v}\in W^{1,s}_0(\Omega)^d,\\
	\inf_{\overline{q}\in M^n}\|q-\overline{q}\|_{L^{s}(\Omega)} &\rightarrow 0 \quad \text{ as }n\rightarrow\infty\quad \forall\,q\in L^s(\Omega),\\
	\inf_{\overline{\bm{\sigma}}\in \Sigma^n}\|\bm{\sigma}-\overline{\bm{\sigma}}\|_{L^{s}(\Omega)} &\rightarrow 0 \quad \text{ as }n\rightarrow\infty\quad \forall\,\bm{\sigma}\in L^s(\Omega)^{d\times d}.
\end{align*}
\end{Assumption}
\begin{Assumption}[Fortin Projector $\Pi^n_\Sigma$]\label{as:Projector_Stress}
For each $n\in\mathbb{N}$ there is a linear projector $\Pi^n_\Sigma \colon L_{\sym}^1(\Omega)^{d\times d}\to \Sigma^n$ such that:
\begin{itemize}
	\item (Preservation of divergence). For any $\bm{\sigma}\in \Lsym{1}$ we have that
\begin{equation*}
\int_\Omega \bm{\sigma}:\BD(\bv) = \int_\Omega \Pi^n_\Sigma(\bm{\sigma}):\BD(\bm{v}) \quad \forall\, \bm{v}\in \Vndiv;
\end{equation*}
\item ($L^{s}$--stability). For every $s\in (1,\infty)$ there is a constant $c>0$, independent of $n$, such that:
\begin{equation*}
	\|\Pi^n_\Sigma \bm{\sigma}\|_{L^{s}(\Omega)} \leq c \|\bm{\sigma}\|_{L^{s}(\Omega)}\qquad \forall\, \bm{\sigma}\in \Lsym{s}.
\end{equation*}
\end{itemize}
\end{Assumption}
\begin{Assumption}[Fortin Projector $\Pi^n_V$]\label{as:Projector_Velocity}
For each $n\in\mathbb{N}$ there is a linear projector $\Pi^n_V :W^{1,1}_0(\Omega)^d\rightarrow V^n$ such that the following properties hold:
\begin{itemize}
\item (Preservation of divergence). For any $\bm{v}\in W^{1,1}_0(\Omega)^d$ we have that
\begin{equation*}
\int_\Omega q\,\diver\bv = \int_\Omega q\,\diver(\Pi^n_V\bv) \quad\, \forall\, q\in M^n.
\end{equation*}
\item ($W^{1,s}$--stability). For every $s\in (1,\infty)$ there is a constant $c>0$, independent of $n$, such that:
\begin{equation*}
\|\Pi^n_V\bm{v}\|_{W^{1,s}(\Omega)} \leq c \|\bm{v}\|_{W^{1,s}(\Omega)}\qquad \forall \,\bm{v}\in W^{1,s}_0(\Omega)^d.
\end{equation*}
\end{itemize}
\end{Assumption}
\begin{Assumption}[Projector $\Pi^n_M$]\label{as:Projector_Pressure}
	For each $n\in\mathbb{N}$ there is a linear projector $\Pi^n_M :L^1(\Omega)\rightarrow M^n$ such that for all $s\in (1,\infty)$ there is a constant $c>0$, independent of $n$, such that:
\begin{align*}
	\|\Pi^n_M q\|_{L^{s}(\Omega)} &\leq c \|q\|_{L^{s}(\Omega)} && \forall\, q\in L^s(\Omega).
\end{align*}
\end{Assumption}

An important consequence of Assumptions \ref{as:Projector_Stress} and \ref{as:Projector_Velocity} is inf-sup stability; i.e.\ we have that for any $s\in(1,\infty)$ there exist two constants $\beta_s,\gamma_s>0$, that do not depend on $n$, such that
\begin{subequations}\label{eq:infsup}
\begin{align}
	\adjustlimits \inf_{q\in M^n\setminus\{0\}}\sup_{\bv\in V^n\setminus\{0\}} \frac{\int_\Omega q\,\diver\bv}{\|\bv\|_{W^{1,s}(\Omega)}\|q\|_{L^{s'}(\Omega)}} \geq \beta_s,\label{eq:infsupvel}\\
	\adjustlimits\inf_{\bv\in V^n_{\diver}\setminus\{0\}}\sup_{\Btau\in \Sigma^n_{\sym}\setminus\{0\}}
\frac{\int_\Omega \bm{\tau}\fp \BD(\bm{v})}{\|\bm{\tau}\|_{L^{s'}(\Omega)}\|\bm{v}\|_{W^{1,s}(\Omega)}} \geq \gamma_s.\label{eq:infsupstress}
\end{align}
\end{subequations}

There are several families of velocity and pressure pairs $V^n$--$M^n$ that satisfy the assumptions above. We can for example mention the Taylor--Hood $\mathbb{P}_k$--$\mathbb{P}_{k-1}$ and MINI elements (see e.g.\ \cite{Boffi2013,Girault1986}). The Scott--Vogelius element $\mathbb{P}_k$--$\mathbb{P}_{k-1}^{\mathrm{disc}}$ is also known to satisfy the assumptions, for instance, on barycentrically refined meshes \cite{Qin1994}; this element has in addition the useful property that discretely divergence-free velocities are in fact also pointwise divergence-free. Regarding the shear stress, if the velocity element consists of globally continuous piecewise polynomials of degree $k\in\mathbb{N}$, then an example of a space satisfying Assumption \ref{as:Projector_Stress} is given by (see e.g.\ \cite{Farrell2019}):
\begin{equation}\label{eq:space_discrete_stresses}
	\Sigma^n = \{\bm{\sigma}\in \Lsym{\infty} \, :\, \bm{\sigma}|_K\in\mathbb{P}_{k-1}(K)^{d\times d},\text{ for all }K\in \mathcal{T}_n\}.
\end{equation}

In the finite element formulation of the (regularised) problem we then look for functions $(\BS^n,\bu^n,p^n)\in \Signsym \times V^n \times M^n_0$ such that:
\begin{subequations}\label{eq:FEFormulation}
\begin{alignat}{2}
	\int_\Omega &\BG_\varepsilon(\BS^n,\BD(\bu^n))\fp  \Btau =0  &\forall\,\Btau \in \Sigma^n,\label{eq:discrete_CR}\\
\alpha\int_\Omega \bu^n\cdot \bv +
	\int_\Omega &\BS^n\fp\BD(\bv) -\int_\Omega p^n  \diver\bv = \int_\Omega \bm{f}\cdot\bv \quad & \forall\,\bv\in V^n, \label{eq:discrete_momentum}\\
		    &-\int_\Omega q  \diver\bu^n = 0 & \forall\, q\in M^n.\label{eq:discrete_mass}
\end{alignat}
\end{subequations}
The Assumptions \ref{as:approximability}--\ref{as:Projector_Pressure} and (G1)--(G4) guarantee that the discrete problem \eqref{eq:FEFormulation} has a solution, and that, as $n\to \infty$, the sequence of discrete solutions converges (weakly) to a solution of the continuous problem (see \cite[Rem.\ 3.8]{Farrell2019}). The rest of the paper will focus on how to solve the nonlinear problems described in this section by means of a semismooth Newton method.

\section{Semismooth Newton method}\label{sec:semismoothNewton}
\subsection{Semismooth functions}
Semismoothness is a concept originally introduced in \cite{Mifflin1977} for the analysis of finite-dimensional nonsmooth optimisation problems (see also \cite{Qi1993,Qi1993a}); it was helpful in the development of a generalised Newton method for functions that are not Fr\'{e}chet-differentiable.

Suppose $F\colon \R^m \to \R^k$ is a locally Lipschitz function. Since $F$ is (by Rademacher's theorem) differentiable everywhere except maybe on a set $U_F$ of zero measure, one can define a generalised notion for the Jacobian of $F$ at any point $z\in \R^m$, called Clarke's differential, as
\begin{equation}\label{eq:ClarkeJac}
	\partial F(z) := \co\{M\in \R^{k\times m} \colon \exists \{z_i\}_i \subset \R^m \setminus U_F \textrm{ with }z_i\to z,\, \nabla F(z_i)\to M \},
\end{equation}
where $\co A$ denotes the convex hull of the set $A$. It is known that Clarke's differential $\partial F(z)$ is nonempty, compact and convex, and that when $F$ is continuously differentiable then $\partial F(z) = \{\nabla F(z)\}$ (for more details the reader is referred to \cite{Clarke1983}). We then say that $F$ is \emph{semismooth at }$z$ if it is directionally differentiable and
\begin{equation}\label{eq:semismoothness}
	\max_{M\in \partial F (z+h)}|f(z+h) -F(z) -Mh| = \mathrm{o}(|h|) \quad \text{as }h\to 0.
\end{equation}
An important consequence of this definition is that if $F$ is a locally Lipschitz and semismooth function in a neighborhood around $z_*\in \R^m$, where $F(z_*)=0$, and $\partial F(z_*)$ is nonsingular, then the generalised Newton iteration
\begin{equation}\label{eq:genNewton}
	z^{j+1} = z^j - M_j^{-1}F(z^j),
\end{equation}
where $M_j\in \partial F(z^j)$ is arbitrary, converges superlinearly to $z_*$, provided $|z^0 - z_*|$ is sufficiently small \cite{Chen2000}.

At the infinite dimensional level, where $F\colon Z \to X$ is now a function defined between two Banach spaces $Z$ and $X$, Rademacher's theorem does not apply and therefore a characterisation of a generalised gradient such as \eqref{eq:ClarkeJac} is not available. Thus, one has to work with a generally weaker and more abstract notion of generalised Jacobians. Let $U$ be a neighborhood of $z\in Z$, and suppose $\partial F\colon U\subset Z \rightrightarrows \mathcal{L}(Z;X)$ is a set-valued function with $\partial F(\hat{z}) \neq \emptyset$ for every $\hat{z}\in U$. We say that $F$ is \emph{$\partial F$-semismooth at }$z$ (in the sense of Ulbrich \cite{Ulbrich2011}) if
\begin{equation}\label{eq:semismoothness_inf}
	\sup_{M\in\partial F(z + h)}\|F(z+h) - F(z) - Mh\|_X = \mathrm{o}(\|h\|_Z)
	\quad \text{as }h\to 0.
\end{equation}
Similarly to the finite-dimensional case, one has that if $F$ is continuously differentiable then it is $\{F'\}$-semismooth. Furthermore, as stated in the following proposition, the semismooth Newton algorithm also leads to superlinear convergence.

\begin{prpstn}[{\cite[Prop.\ 2.7]{Ulbrich2003}}]\label{prop:Newton_convergence}
	Let $X$, $Z$ be Banach spaces, and let $U$ be a neighborhood of a point $z_*\in Z$. Suppose that the function $F\colon Z\to X$ is $\partial F$-semismooth at $z_*$ and that $F(z_*)=0$. Suppose that there exists a positive constant $c$, such that for every $z\in U$ and every $M\in \partial F(z)$, $M$ is invertible with $\|M^{-1}\|_{\mathcal{L}(X;Z)}\leq c$. Then the generalised Newton--Kantorovich iteration defined by \eqref{eq:genNewton} converges superlinearly to $z_*$ for all $z^0$ sufficiently close to $z_*$.
\end{prpstn}

We will now formulate a semismooth Newton method for the problem \eqref{eq:PDE_}. Define for $r\in (1,\infty)$ the mapping
\begin{equation}\label{eq:def_F}
\begin{gathered}
	F\colon Z:= \Lsymtr{r'}\times V^r \times L^{r'}_0(\Omega)\to
X:=	\Lsymtr{q} \times W^{-1,r'}(\Omega)^d  \times L_0^r(\Omega),\\
	F(\BS,\bu,p) := \begin{pmatrix}
		\overline{\BG}(\BS,\BD(\bu))\\
		\alpha \bu  -\diver \BS  + \nabla p - \bm{f}\\
		\diver \bu
	\end{pmatrix},
\end{gathered}
\end{equation}
where $q>1$  and $\overline{\BG}$ denotes the Nemytskii operator associated to $\BG$. Observe that the second and third entries in the definition of $F$ are given by linear and bounded operators, and so the semismoothness (or lack thereof) of $F$ depends solely on that of $\overline{\BG}$. This is advantageous because from this one immediately obtains a candidate for the multivalued function $\partial F$ from Proposition \ref{prop:Newton_convergence}. Namely, for $(\BS,\bu,p)\in Z$, the generalised Jacobian $\partial F(\BS,\bu,p)$ will consist of the linear operators $M\in \mathcal{L}(Z;X)$ of the form
\begin{equation}\label{eq:gen_Jacobian_F}
	M[\BT,\bv,q] =
	\begin{pmatrix}
		\bm{d}_1 \BT + \bm{d}_2\BD(\bv)\\
		\alpha \bv - \diver \BT + \nabla \txtbl{m}\\
		\diver \bv
	\end{pmatrix},
\end{equation}
where $[\bm{d}_1,\bm{d}_2]$ is a measurable selection of $\partial \BG(\BS(\cdot),\BD(\bu(\cdot)))$, with $\partial \BG$ being the Clarke generalised gradient of $\BG$ (cf.\ \cite[Def.\ 3.40]{Ulbrich2011}). The main algorithm can be summarised as follows:
\begin{lgrthm}\label{alg:algorithm}
	Select an initial guess $(\BS_0,\bu_0,p_0)\in \Lsym{r'}\times V^r \times \Lmean{r'}$, and for each $k=0,1,2,\ldots$ perform the steps:
	\begin{enumerate}
		\item If $F(\BS_k,\bu_k,p_k)=0$, then terminate with $(\BS_k,\bu_k,p_k)$;
		\item Choose an arbitrary element $M_k\in \partial F$, where $\partial F$ is defined through \eqref{eq:gen_Jacobian_F};
		\item Solve the system $M_k[\tilde{\BS}_k,\tilde{\bu}_k,\tilde{p}_k] = -F(\BS_k,\bu_k,p_k)$;
		\item Update $(\BS_k,\bu_k,p_k) \gets (\BS_k,\bu_k,p_k) + (\tilde{\BS}_k,\tilde{\bu}_k,\tilde{p}_k)$.
	\end{enumerate}
\end{lgrthm}

The continuity of Nemytskii mappings between Lebesgue spaces has been well studied  and it is known (see e.g.\ \cite[Thm. 4]{Goldberg1992}) that if $\BG$ satisfies, for some $q>1,$ the growth condition
\begin{equation}\label{eq:growth_G}
	|\BG(\Bsigma,\Btau)| \leq c(|\Bsigma|^{r'/q} +
	|\Btau|^{r/q}) + \tilde{c}
	\qquad \forall\, \Bsigma,\Btau\in \Rds,
\end{equation}
where $c,\tilde{c}>0$, then the Nemytskii map
\begin{equation}\label{eq:Nemytskii}
	(\BS,\BT) \in \Lsym{r'}\times \Lsym{r} \mapsto
	\overline{\BG}(\BS,\BT):=	\BG(\BS(\cdot),\BT(\cdot)) \in \Lsym{q},
\end{equation}
is continuous. Suppose for a moment that the partial derivatives
\begin{equation*}
	\partial_{\Bsigma}\BG,\partial_{\Btau}\BG\colon \Rds\times \Rds\to \mathcal{L}(\Rds;\Rds),
\end{equation*}
are continuous, and define
\begin{equation}\label{eq:def_q}
	q := \min\{r,r'\} - \varepsilon \in [1,\min\{r,r'\}),\text{  with }\varepsilon\in (0,\min\{r,r'\}-1].
\end{equation}
If we assume that the following growth conditions hold:
\begin{equation}\label{eq:growth_derivatives}
	\begin{split}
		\|\partial_{\Bsigma}\BG(\Bsigma,\Btau)\|_{\mathcal{L}(\Rdd;\Rdd)} &\leq c_1(|\Bsigma|^{\frac{r'}{s}} + |\Btau|^{\frac{r}{s}}) + c_2,\\
		\|\partial_{\Btau}\BG(\Bsigma,\Btau)\|_{\mathcal{L}(\Rdd;\Rdd)} &\leq \tilde{c}_1(|\Bsigma|^{\frac{r'}{p}} + |\Btau|^{\frac{r}{p}}) + \tilde{c}_2,
	\end{split}
\end{equation}
where $c_i,\tilde{c}_i$ are non-negative constants, and
\begin{equation}\label{eq:def_s_p}
	s := \frac{qr'}{r'-q}\in [r,\infty)\qquad p:=\frac{qr}{r-q}\in \left[\frac{r}{r-1},\frac{r}{r-2} \right),
\end{equation}
then the operators
\begin{equation*}
	(\BS,\BT)\in \Lsym{r'}\times \Lsym{r} \mapsto \partial_{\Bsigma}\BG(\BS(\cdot),\BT(\cdot)),\partial_{\Btau}\BG(\BS(\cdot),\BT(\cdot))\in \Lsym{q},
\end{equation*}
are continuous, which in turn implies that the corresponding Nemytskii maps are continuously differentiable (see e.g.\ \cite[Thm. 7]{Goldberg1992}), and the derivatives
\begin{equation}
	\begin{split}
		\partial_{1}\overline{\BG}\colon \Lsym{r'}\times \Lsym{r} &\mapsto \mathcal{L}(\Lsym{r'};\Lsym{q}),\\
		\partial_{2}\overline{\BG}\colon \Lsym{r'}\times \Lsym{r} &\mapsto \mathcal{L}(\Lsym{r};\Lsym{q}),
\end{split}
\end{equation}
are given by
\begin{equation}
	\begin{split}
		(\partial_{1}\overline{\BG} (\BS,\BT)\tilde{\BS})(\cdot) &= \partial_{\Bsigma}\BG(\BS(\cdot),\BT(\cdot))\tilde{\BS}(\cdot), \quad \tilde{\BS} \in \Lsym{r'},\\
		(\partial_{2}\overline{\BG} (\BS,\BT)\tilde{\BT})(\cdot) &= \partial_{\Btau}\BG(\BS(\cdot),\BT(\cdot))\tilde{\BT}(\cdot), \quad \tilde{\BT} \in \Lsym{r}.\\
\end{split}
\end{equation}
In summary, if one assumes that the partial derivatives of $\BG$ are continuous, and they satisfy the growth assumptions \eqref{eq:growth_derivatives}, then the Nemytskii map \eqref{eq:Nemytskii} is Fr\'{e}chet-differentiable (and therefore so is $F$), which would make the use of a Newton--Kantorovich algorithm possible. Note that the number $q$ defined \eqref{eq:def_q} is smaller than both $r$ and $r'$; this is an example of a norm gap, which is a phenomenon that appears in other contexts such as optimal control theory (see e.g.\ \cite{Troeltzsch2010}).

As mentioned in the introduction, we wish here to work with Bingham-like constitutive relations that do not satisfy the differentiability requirement just described, and we therefore wish to have at our disposal a similar characterisation for semismoothness of superposition operators. 
The first results in this direction were obtained in \cite{Ulbrich2003} assuming Lipschitz continuity, which is too restrictive for our setting. Thankfully, this assumption was relaxed in \cite{Schiela2008}, leading to conditions that guarantee semismoothness that are more widely applicable. The following lemma, consequence of \cite[Prop. A.1]{Schiela2008}, summarises the required assumptions in the current setting.

\begin{lmm}\label{lm:lemma_semismooth}
	Let $(\BS,\BT)\in \Lsym{r'}\times \Lsym{r}$, where $r\in (1,\infty)$ is arbitrary, and let $q\geq 1$ be defined by \eqref{eq:def_q}. Suppose that $\BG\colon \Rds\times \Rds \to \Rds$ is a locally Lipschitz function satisfying the growth assumption \eqref{eq:growth_G}, with partial derivatives $\partial_{\Bsigma}\BG$, $\partial_{\Btau}\BG$ satisfying \eqref{eq:growth_derivatives}. Suppose that $M:=[\bm{d}_1 \,\,\, \bm{d}_2]\colon \Rds\times \Rds \to \mathcal{L}(\Rds\times \Rds;\Rds)$ satisfies
	\begin{equation}\label{eq:growth_lemma}
	\begin{split}
		\|\bm{d}_1(\Bsigma,\Btau)\|_{\mathcal{L}(\Rds;\Rds)} &\leq c_1(|\Bsigma|^{\frac{r'}{s}} + |\Btau|^{\frac{r}{s}}) + c_2,\\
		\|\bm{d}_2(\Bsigma,\Btau)\|_{\mathcal{L}(\Rds;\Rds)} &\leq c_1(|\Bsigma|^{\frac{r'}{p}} + |\Btau|^{\frac{r}{p}}) + c_2,
	\end{split}
\end{equation}
where $c_1,c_2$ are positive constants, and $s,p$ are defined by \eqref{eq:def_s_p}; assume further that $M$ can be written as the a.e.\ pointwise limit of continuous functions. Then, if $\BG$ is semismooth, the associated Nemytskii operator $\overline{\BG}$ defined by \eqref{eq:Nemytskii} will satisfy
\begin{equation}\label{eq:lemma_semismooth}
	\|\overline{\BG}(\BS + \tilde{\BS},\BT + \tilde{\BT}) -
	\overline{\BG}(\BS,\BT)
	- M(\BS+\tilde{\BS},\BT+\tilde{\BT})(\tilde{\BS},\tilde{\BT})\|_{L^q(\Omega)}
	= \mathrm{o}(\|(\tilde{\BS},\tilde{\BT})\|_{L^{r'}(\Omega)\times L^{r}(\Omega)})
	\quad \text{as }(\tilde{\BS},\tilde{\BT})\to \bm{0}.
\end{equation}
\end{lmm}
\begin{proof}
	First note that the existence of the partial derivatives $\partial_{\Bsigma}\BG$, $\partial_{\Btau}\BG$ is guaranteed by Rademacher's theorem. Define now  the local Lipschitz constants of $\BG$ as
	\begin{equation}
		L_1(\Bsigma,\Btau):= \sup_{|\Bsigma - \tilde{\Bsigma}|\leq 1}
		\frac{|\BG(\Bsigma,\Btau) - \BG(\tilde{\Bsigma},\Btau)|}{|\Bsigma - \tilde{\Bsigma}|},\qquad
L_2(\Bsigma,\Btau):= \sup_{|\Btau - \tilde{\Btau}|\leq 1}
\frac{|\BG(\Bsigma,\Btau) - \BG(\Bsigma,\tilde{\Btau})|}{|\Btau - \tilde{\Btau}|},
	\end{equation}
	where $\Bsigma,\tilde{\Bsigma},\Btau,\tilde{\Btau}\in \Rds$; then the growth assumption \eqref{eq:growth_derivatives} can be reframed as
\begin{equation}
	\begin{split}
		|L_1(\Bsigma,\Btau)| &\leq c_1(|\Bsigma|^{\frac{r'}{s}} + |\Btau|^{\frac{r}{s}}) + c_2,\\
		|L_2(\Bsigma,\Btau)| &\leq c_1(|\Bsigma|^{\frac{r'}{p}} + |\Btau|^{\frac{r}{p}}) + c_2.\\
	\end{split}
\end{equation}
Since $M$ is a pointwise limit of continuous functions (i.e.\ it is a Baire--Carath\'{e}odory function, but without an explicit space-dependence),  \cite[Prop.\ A.1]{Schiela2008} implies that the function
\begin{equation}
  \renewcommand{\arraystretch}{1.3}
\left\{
	\begin{array}{ccc}
		\frac{M(\BS(x)+\Bsigma,\BT(x)+\Btau)[\Bsigma,\Btau] - (\BG(\BS(x)+\Bsigma,\BT(x)+\Btau) - \BG(\BS(x),\BT(x))) }{|(\Bsigma,\Btau)|} & \Longleftrightarrow & (\Bsigma,\Btau)\neq \bm{0}, \\
	\bm{0} & \Longleftrightarrow & (\Bsigma,\Btau)=\bm{0},
	\end{array}
\right.
\end{equation}
 defined for $(x,\Bsigma,\Btau)\in \Omega \times \Rds\times \Rds$, induces a continuous Nemytskii operator, which in particular implies the desired condition \eqref{eq:lemma_semismooth}.
	\end{proof}

	In practice, the function $M$ in Lemma \ref{lm:lemma_semismooth} arises from derivatives of $\BG$, and so the growth conditions \eqref{eq:growth_derivatives} and \eqref{eq:growth_lemma} are redundant. This is the case for instance with the generalised Jacobian defined by \eqref{eq:gen_Jacobian_F}. After having identified the necessary assumptions to guarantee semismoothness, it only remains to examine the solvability of the linearised system.
	\begin{Assumption}\label{as:linearised}
		Let $(\BS,\bu)\in \Lsym{r'}\times V^r$, with $r\in (1,\infty)$, and define $q$ through \eqref{eq:def_q}. Let $\BG\colon \Rds\times \Rds\to \Rds$ be locally Lipschitz. Given a triplet $(\BH,\bh,h)\in \Lsym{q}\times L^{r'}(\Omega)^d\times \Lmean{r}$, there exists a solution $(\BT,\bv,\txtbl{m})\in \Lsym{r'}\times V^r\times \Lmean{r'}$ of the linear system
\begin{align}
	\bm{d}_1 \BT + \bm{d}_2& \BD(\bv) = \BH, \notag \\
	\alpha \bv - \diver \BT &+ \nabla \txtbl{m} = \bh,\label{eq:linearised_PDE}\\
	\diver \bv &= h, \notag
\end{align}
where $[\bm{d}_1,\bm{d}_2]$ is a measurable selection of $\partial \BG(\BS(\cdot),\BD\bu(\cdot))$, such that
\begin{equation}\label{eq:estimate_linearised}
	\|\BT\|_{L^{r'}(\Omega)}
	+ \|\bv\|_{W^{1,r}(\Omega)}
	+ \|\txtbl{m}\|_{L^{r'}(\Omega)}
	\leq c(\|\BH\|_{L^q(\Omega)}
	+ \|\bh\|_{L^{r'}(\Omega)}
	+ \|h\|_{L^r(\Omega)}),
\end{equation}
where $c$ is a positive constant.
	\end{Assumption}
	An immediate consequence of Proposition \ref{prop:Newton_convergence} is the following convergence result.
\begin{thrm}
	Let $\BG\colon \Rds\times \Rds \to \Rds$ be a function satisfying assumptions (G1)--(G4) for some $r>1$. Suppose that the growth conditions \eqref{eq:growth_G} and \eqref{eq:growth_derivatives} hold, and that Assumption \ref{as:linearised} is satisfied. Let $(\BS,\bu,p)\in \Lsym{r'}\times V^r\times \Lmean{r'}$ be such that $F(\BS,\bu,p)=0$, where $F$ is defined by \eqref{eq:def_F}. Then the iterates produced by Algorithm \ref{alg:algorithm} converge superlinearly to $(\BS,\bu,p)$, for any initial guess $(\BS_0,\bu_0,p_0)$ sufficiently close to $(\BS,\bu,p)$. 
\end{thrm}

	In some cases higher regularity results allow for less restrictive growth conditions for $\BG$, meaning that a larger $q$ could be selected in \eqref{eq:def_F}. 
For instance, the strict monotonicity \eqref{eq:coercivity_monotonicity_regularised} might result in higher regularity for the problem associated to the regularisation $\BG_\varepsilon$ defined by \eqref{eq:CR_regularised}, namely $(\BS_\varepsilon,\bu_\varepsilon)\in W_{\mathrm{sym}}^{1,\min\{2,r'\}}(\Omega)^{d\times d}\times W^{2,\min\{2,r\}}(\Omega)^d$; in that case one could e.g.\ guarantee the continuity of the Nemytskii operator $\overline{\BG}_\varepsilon$ with the more general growth condition
	\begin{equation}\label{eq:growth_higher_reg}
		|\BG_\varepsilon(\Bsigma,\Btau)| \leq c\left(
		|\Bsigma|^{\frac{\min\{2,r'\}^*}{q}}
		+ |\Btau|^{\frac{\min\{2,r\}^*}{q}}\right) + \tilde{c}
		\qquad c,\tilde{c}>0,
		\qquad \Bsigma,\Btau\in \Rds,
	\end{equation}
	where  $p^*:= \frac{p d}{d-p}$ denotes the Sobolev exponent of $p$ (it is defined as infinity when $p=d$). 

	Without assuming any sort of strict monotonicity of the constitutive relation, it might prove difficult to solve the linearised system \eqref{eq:linearised_PDE} from Assumption \ref{as:linearised}. In this case the use of the regularisation \eqref{eq:CR_regularised} also brings an advantage: the monotonicity condition \eqref{eq:coercivity_monotonicity_regularised} implies by the implicit function theorem that the partial derivatives of $\BG_\varepsilon$ are uniformly positive definite, and so the matrix function $A_\varepsilon := -(\bm{d}_1^\varepsilon)^{-1}\bm{d}_2^\varepsilon$ is uniformly positive definite, which means that the system \eqref{eq:linearised_PDE} can be rewritten as a Stokes-type system:
	\begin{equation}\label{eq:linearised_stokes}
\begin{aligned}
	\alpha\bv_\varepsilon - \diver(A_\varepsilon \BD(\bv_\varepsilon)) &+ \nabla \txtbl{m}_\varepsilon
= \bh + \diver \tilde{\BH},\\
\diver \bv_\varepsilon &= h,
\end{aligned}
\end{equation}
where $\tilde{\BH} := (\bm{d}^\varepsilon_1)^{-1}\BH \in L^{q}(\Omega)^{d\times d}$. If $q\geq \max\{r,r'\}$, the system \eqref{eq:linearised_stokes} can be solved by standard energy arguments, and the required estimate \eqref{eq:estimate_linearised} is for instance  a consequence of \cite[Thm. 1.4]{Bulicek2016a}. In the opposite case $q\in (1,\max\{r,r'\})$ the situation becomes slightly more delicate since the term $\diver \tilde{\BH}$, with $\BH\in \Lsym{q}$ and $q$ small, does not allow the use of standard energy estimates. The results from \cite{Bulicek2016a} still guarantee the existence of a solution, but the estimate takes the form
\begin{equation}
	\|\BT_\varepsilon\|_{L^q(\Omega)} + \|\bv_\varepsilon\|_{W^{1,q}(\Omega)} + \|\txtbl{m}_\varepsilon\|_{L^{q}(\Omega)} \leq c(1+ \|\BH\|_{L^q(\Omega)} + \|\bh\|_{L^{r'}(\Omega)}
	+\|h\|_{L^r(\Omega)}).
\end{equation}
In this case, the addition of a smoothing step to Algorithm \ref{alg:algorithm} that takes the solution back to $(\BT_\varepsilon,\bv_\varepsilon,\txtbl{m}_\varepsilon)\in \Lsym{r'}\times W^{1,r}_0(\Omega)^d\times \Lmean{r'}$ would be necessary (see e.g.\ \cite{Ulbrich2011} for more details on smoothing operators).

\begin{crllr}\label{cor:corollary}
	Let $\BG\colon \Rds\times \Rds \to \Rds$ be a function satisfying (G1)--(G4) for some $r>1$, and let $(\BS_\varepsilon,\bu_\varepsilon,p_\varepsilon)\in \Lsym{r}\times V^r\times \Lmean{r'}$ be such that $F_\varepsilon(\BS_\varepsilon,\bu_\varepsilon,p_\varepsilon)=0$, where $F_\varepsilon$ is defined through \eqref{eq:def_F} with $\BG$ replaced by $\BG_\varepsilon$. Then, if $q\geq \max\{r,r'\}$ (recall \eqref{eq:growth_higher_reg}), the iterates produced by Algorithm \ref{alg:algorithm} (employing $F_\varepsilon$ instead of $F$) converge superlinearly to $(\BS_\varepsilon,\bu_\varepsilon,p_\varepsilon)$, for any initial guess $(\BS_0^{\varepsilon},\bu_0^\varepsilon,p_0^\varepsilon)$ sufficiently close to $(\BS_\varepsilon,\bu_\varepsilon,p_\varepsilon)$. If $q<\max\{r,r'\}$ the same result holds if a smoothing step is included in Algorithm \ref{alg:algorithm}. 
\end{crllr}

\begin{rmrk}
	While the previous discussion dealt with the continuous problem \eqref{eq:PDE}, the exact same arguments regarding e.g.\ continuity could be applied to a weakly enforced form of the constitutive relation
	\begin{equation*}
		\int_\Omega \BG(\BS,\BD(\bu))\fp \Btau = \bm{0} \qquad
		\forall \, \Btau \in C_0^\infty(\Omega)^{d\times d},
	\end{equation*}
	and so in turn also to the finite element discretisation \eqref{eq:FEFormulation}. At the discrete level all norms are equivalent, which would make it straightforward to obtain estimate \eqref{eq:estimate_linearised}, without needing sophisticated results like the ones from \cite{Bulicek2016a}. However, it is likely that this would then result in mesh-dependent behaviour of the iterations; it is for this reason that we endeavoured to carry out the analysis at the function space level.
\end{rmrk}

\begin{rmrk}
The convective term can be incorporated into the analysis in a similar manner as is done for the Newtonian problem. For example, if we employ a Picard-like linearisation of the form
\begin{equation*}
	\mathcal{B}[\bu;\bv,\bm{\phi}]	:= -\int_\Omega \bv\otimes \bu\fp \bm{\phi},
\end{equation*}
where $\bu\in V^r$ is the current solution and $\bm{\phi}\in C^\infty_0(\Omega)^d$ is a test function, the solvability of the linearised system remains intact because of the skew-symmetry property $\mathcal{B}[\bu;\bv,\bv]=0$ for all $\bv\in V^r$ (but note that one has to assume that $r$ is large enough so that choosing $\bm{\phi}=\bv$ is allowed), which is a consequence of the divergence-free condition $\diver \bu = 0$. If one wishes to employ instead the form of the convective term that would arise from Newton linearisation:
\begin{equation*}
	\mathcal{B}[\bu;\bv,\bm{\phi}]	:= -\int_\Omega (\bv\otimes \bu + \bu\otimes \bv)\fp \bm{\phi},
\end{equation*}
the situation is slightly more delicate, and even in the Newtonian case it is necessary to assume e.g.\ that the viscosity is large enough in order to guarantee that the linearised system can be solved (see e.g.\ \cite[Lm.\ 3.2]{Girault1986}).
\end{rmrk}

\subsection{The Bingham constitutive relation}\label{sec:Bingham}
In this section we will look at some explicit examples of constitutive relations that fall into the framework described in the previous section. Namely, we will show that the Bingham constitutive relation can be described by means of a semismooth function $\BG$ that also satisfies appropriate growth conditions.

First, we claim that the relation \eqref{eq:Bingham_implicit_2} defining a Bingham fluid is semismooth. Note that, excepting the first term, all the terms appearing in the following function
\begin{equation}\label{eq:Bingham_implicit_2a}
\BG\colon \Rds\times \Rds \to \Rds \qquad
\BG(\Bsigma,\Btau) := |\Btau|\Bsigma - (\tau_* + 2\nu |\Btau|)\Btau\qquad \tau_*,\nu>0,
\end{equation}
are smooth. Therefore we only need concern ourselves with the function
\begin{equation}\label{eq:h_semismooth}
	\bh\colon \Rds\times\Rds \to \Rds\qquad
	\bh(\Bsigma,\Btau) := |\Btau|\Bsigma.
\end{equation}
The function $\bh$ is locally Lipschitz, and its (Clarke) generalised Jacobian is given by
\begin{equation}\label{eq:C_jacobian_h}
\renewcommand{\arraystretch}{1.5}
\partial \bh(\Bsigma,\Btau)=
\left\{
\begin{array}{cc}
	\left\{\left[|\Btau|\BI, \Bsigma\otimes \frac{\Btau}{|\Btau|}\right]\right\}, & \textrm{ if } |\Btau|\neq 0,\\
	\left\{\left[ \bm{0} , \Bsigma\otimes\Bphi\right]\, :\, \Bphi\in \Rds,\, |\Bphi|\leq 1  \right\}	  , & \textrm{ if }|\Btau|=0,\\
\end{array}
\right.
\end{equation}
where $\BI$ denotes the fourth order identity tensor. Since for $\Btau \neq \bm{0}$ the function $\bh$ is actually smooth, we only need to prove semismoothness at $\Btau = \bm{0}$:
\begin{equation*}
|\bh(\Bsigma + \tilde{\Bsigma},\tilde{\Btau}) - \bh(\Bsigma,\bm{0})
- |\tilde{\Btau}|\tilde{\Bsigma} - |\tilde{\Btau}|(\Bsigma +\tilde{\Bsigma}) |
= |\tilde{\Btau}||\tilde{\Bsigma}||(\tilde{\Bsigma},\tilde{\Btau})| = \textrm{o}(|(\tilde{\Bsigma},\tilde{\Btau})|),
	\text{ as }|(\tilde{\Bsigma},\tilde{\Btau})|\to 0.
\end{equation*}
The function \eqref{eq:Bingham_implicit_2a} is therefore semismooth and its generalised Jacobian is given by
\begin{equation}\label{eq:C_jacobian_G}
\renewcommand{\arraystretch}{1.5}
\partial \BG(\Bsigma,\Btau)=
\left\{
\begin{array}{cc}
	\left\{\left[|\Btau|\BI , \Bsigma\otimes \frac{\Btau}{|\Btau|} - (\tau_* + |\Btau|)\BI - \frac{\Btau\otimes \Btau}{|\Btau|} \right]\right\}, & \textrm{ if } |\Btau|\neq 0,\\
	\left\{\left[ \bm{0} , \Bsigma\otimes\Bphi -\tau_*\BI \right]\, :\, \Bphi\in \Rds,\, |\Bphi|\leq 1  \right\}	  , & \textrm{ if }|\Btau|=0.\\
\end{array}
\right.
\end{equation}
For the Bingham constitutive relation we have $r=2$ and so the solution to the problem belongs to $(\BS,\bu,p)\in \Lsym{2}\times W_0^{1,2}(\Omega)^d\times \Lmean{2}$. On the other hand, a simple application of Young's inequality yields the growth conditions for $\BG$ and its derivatives:
\begin{align}
	|\BG(\Bsigma,\Btau)| &\leq c(|\Bsigma|^2 + |\Btau|^2) + \tilde{c}, \notag\\
	|\partial_{\Bsigma}\BG(\Bsigma,\Btau)| &\leq c(|\Bsigma| + |\Btau|) + \tilde{c}, \label{eq:growth_Bingham_1}\\
	|\partial_{\Btau}\BG (\Bsigma,\Btau)| &\leq c(|\Bsigma|+|\Btau|) + \tilde{c},\notag
\end{align}
where $c,\tilde{c}>0$. The growth conditions \eqref{eq:growth_Bingham_1} alone do not suffice to be able to apply the results from the previous section, since they only imply that the Nemytskii operator \eqref{eq:Nemytskii} is continuous into $\Lsym{1}$, which brings great difficulties since it is unlikely that one can solve the linearised system \eqref{eq:linearised_PDE} with a uniform bound on the solution, if the linearised system is at all solvable (recall the lack of monotonicity).

Thankfully, at this point we can make use of regularity results that guarantee the boundedness of gradients for Bingham and Herschel--Bulkley models \cite{Fuchs1999,Fuchs2000}. Noting also that the mapping $\BD(\bu_\varepsilon)\mapsto \BS_\varepsilon$ defined by the relation \eqref{eq:CR_regularised} is Lipschitz continuous \cite[Lm.\ 4.5]{Maringova2018}, we obtain also boundedness of the stresses, which implies that the corresponding Nemytskii operator $\overline{\BG}_\varepsilon$ will be continuous into $\Lsym{2}$. Moreover, thanks to the uniform monotonicity condition \eqref{eq:coercivity_monotonicity_regularised}, it is possible to solve the linearised system in a straightforward manner, which implies that Corollary \ref{cor:corollary} can be applied to obtain convergence of Algorithm \ref{alg:algorithm} (note that $\BG_\varepsilon$ is semismooth, because it is a composition of a linear transformation with $\BG$).

\begin{rmrk}
	Strictly speaking, the expression \eqref{eq:Bingham_implicit_2a} does not exactly represent the constitutive relation, since when $|\Btau|=0$, the equation $\BG(\Bsigma,\Btau)=\bm{0}$ is satisfied for arbitrary $\Bsigma\in\Rds$, and not just $|\Bsigma|\leq \tau_*$ as required by the Bingham relation \eqref{eq:Bingham_intro}. However, the regularised relation $\BG_\varepsilon$ does not exhibit this unwanted behaviour; e.g.\ Figure \ref{fig:regularised_CR} was obtained using this expression. Numerical experiments in Section \ref{sec:Examples} will also show that the semismooth Newton algorithm will perform well using the regularised version of \eqref{eq:Bingham_implicit_2a}.
\end{rmrk}

The alternative expression \eqref{eq:Bingham_implicit_1} used to define the same Bingham constitutive relation satisfies the same growth conditions \eqref{eq:growth_Bingham_1}, but is not semismooth. To see this, note that its Clarke Jacobian is given by
\begin{equation*}
\renewcommand{\arraystretch}{1.5}
\partial {\BG}(\Bsigma,\Btau)=
\left\{
\begin{array}{cc}
	\left\{\left[(|\Bsigma|-\tau_*)^+\BI
			+ \mathds{1}_{\{|\Bsigma|> \tau_*\}}\frac{(\Bsigma-2\nu\Btau)\otimes \Bsigma}{|\Bsigma|},
			-2\nu(\tau_* + (|\Bsigma| - \tau_*)^+)\BI
	\right]\right\}, & \textrm{ if } |\Bsigma|\neq \tau_*,\\
	\left\{\left[(\Bsigma-2\nu \Btau)\otimes\Bphi, -2\nu\tau_*\BI \right]\, :\, \Bphi\in \Rds,\, |\Bphi|\leq 1  \right\}	  , & \textrm{ if }|\Bsigma|=\tau_*,\\
\end{array}
\right.
\end{equation*}
Take $(\Bsigma,\Btau)\in \Rds\times \Rds$ with $|\Bsigma|=\tau_*$ and $\Btau$ arbitrary. To see that ${\BG}$ is not semismooth, take a sequence $(\Bsigma_k,\Btau_k) \to (\bm{0},\bm{0})$ such that $|\Bsigma + \Bsigma_k| = \tau_*$. Then for any $M_{\Bphi} \in \partial{\BG}(\Bsigma+\Bsigma_k,\Btau + \Btau_k)$ we have that
\begin{equation*}
	|\BG(\Bsigma+\Bsigma_k,\Btau+\Btau_k) - \BG(\Bsigma,\Btau) - M_{\Bphi}[\Bsigma_k,\Btau_k]|
	= |\Bsigma_k\fp \Bphi||2\nu(\Btau + \Btau_k) - (\Bsigma+\Bsigma_k)|,
\end{equation*}
which is not $\mathrm{o}(|(\Bsigma_k,\Btau_k)|)$ as $k\to \infty$. However, it is possible to prove that
\begin{equation*}
|\BG(\Bsigma + \tilde{\Bsigma},\Btau + \tilde{\Btau}) - \BG(\Bsigma,\Btau)
- M[\tilde{\Bsigma},\tilde{\Btau}] |
=  \textrm{o}(|(\tilde{\Bsigma},\tilde{\Btau})|^\alpha),
	\text{ as }|(\tilde{\Bsigma},\tilde{\Btau})|\to 0,
\end{equation*}
for $(\Bsigma,\Btau)\in \Rds\times \Rds$ with $|\Bsigma|=\tau_*$, where $M\in \partial\BG(\Bsigma+\tilde{\Bsigma},\Btau+\tilde{\Btau})$, and $\alpha\in (0,1)$ is arbitrary. This suggests that Algorithm \ref{alg:algorithm} would still be convergent, albeit at a slower rate.

The following is yet another expression that can be used to describe the Bingham relation:
\begin{equation}\label{eq:Bingham_implicit_3}
	\BG(\Bsigma,\Btau) := (|\Bsigma|-\tau_*)^+\frac{\Bsigma}{|\Bsigma|} - 2\nu\Btau.
\end{equation}
The expression \eqref{eq:Bingham_implicit_3} in fact satisfies a less restrictive growth condition
\begin{equation}
	|\BG(\Bsigma,\Btau)|\leq c(|\Bsigma| + |\Btau|) + \tilde{c}
	\qquad c,\tilde{c}>0,
\end{equation}
which e.g. guarantees the continuity of the Nemytskii operator $\overline{\BG}$ into $\Lsym{2}$ without the need for regularity results, but unfortunately it also fails to be semismooth. This shows that different expressions for the same constitutive can relation result in different behaviour, and so it would be worthwhile in future research to look for the most advantageous ways of writing the constitutive relation.

\begin{rmrk}
The Herschel--Bulkley relation can be analysed in a similar manner. For instance, we could use the following expression:
\begin{equation}\label{eq:HB_implicit_2}
\BG\colon \Rds\times \Rds \to \Rds \qquad
\BG(\Bsigma,\Btau) := |\Btau|\Bsigma - (\tau_* + 2\nu |\Btau|^{r-1})\Btau\qquad \tau_*,\nu>0,r>1.
\end{equation}
Since $r>1$, the function $\Btau\mapsto |\Btau|^{r-1}\Btau$ is actually differentiable, and so by the same arguments employed above, the function \eqref{eq:HB_implicit_2} is semismooth. As for the growth condition, one can e.g.\ see that the regularised expression satisfies
\begin{equation}\label{eq:growth_HB}
	|\BG_\varepsilon(\Bsigma,\Btau)| \leq c(|\Bsigma|^{\max\{r,2\}} + |\Btau|^{\max\{r,2\}}) + \tilde{c}\qquad c,\tilde{c}>0.
\end{equation}
Then, for instance if we look at the shear-thinning case $r<2$, the natural higher integrability $(\BS,\Du)\in \Lsym{2^*}\times \Lsym{r^*}$ would suggest we rewrite \eqref{eq:growth_HB} as
\begin{equation}\label{eq:growth_HB2}
	|\BG_\varepsilon(\Bsigma,\Btau)| \leq c(|\Bsigma|^{\frac{2^*}{2^*/2}} + |\Btau|^{\frac{r^*}{r^*/2}}) + \tilde{c}\qquad c,\tilde{c}>0.
\end{equation}
Therefore, if we assume e.g.\ that $r\geq \frac{3d}{d+2}$ (which ensures that $\frac{r^*}{2}\geq r'$), Corollary \ref{cor:corollary} implies the convergence of Algorithm \ref{alg:algorithm} without the need for a smoothing step. We note that the condition $r\geq \frac{3d}{d+2}$ is advantageous when analysing the convergence of finite element approximations for the system including the convective term (see e.g.\ \cite{Diening:2013}).
\end{rmrk}

{\color{black}
\begin{rmrk}
	The method presented in this work has the advantage that it works directly with the basic formulation conventionally derived in continuum mechanics: balance laws in the form of partial differential equations \eqref{eq:PDE_} supplemented with a constitutive relation in the form of a pointwise algebraic constraint \eqref{eq:PDE_CR}. This leads to a lot of flexibility in the choice of constitutive relations; for instance, one could easily capture with the general form \eqref{eq:PDE_CR} constitutive relations of the type
\begin{equation}\label{eq:activated_CR}
  \renewcommand{\arraystretch}{1.2}
\left\{
	\begin{array}{ccc}
		|\BS|\leq \tau_* & \Longleftrightarrow & \Du = \bm{0}, \\
|\BS|> \tau_* & \Longleftrightarrow & \BS = \tau_*\displaystyle\frac{\Du}{|\Du|} + \BS_2,\quad \text{where }\Du = \alpha(|\BS|)\BS.
	\end{array}
\right.
\end{equation}
In contrast, the variational inequality formulation makes heavy use of the form of the potential, e.g.\ for the Bingham model, $ \int_\Omega(\nu |\Du|^2+ \tau_*|\Du|)$ and it is not immediately clear how one should modify the formulation to capture a relation such as \eqref{eq:activated_CR}. In addition, the variational inequality formulation is itself derived as a \emph{consequence} of the laws \eqref{eq:PDE}, and in some cases it is not completely obvious that a solution of the variational inequality is also a solution of the weak form of \eqref{eq:PDE}, e.g.\ in the three-dimensional unsteady Bingham model including inertia, where testing with the velocity $\bu$ is not allowed.
\end{rmrk}
}

\section{Numerical examples}\label{sec:Examples}
In this section we implement Algorithm \ref{alg:algorithm} with the regularisation \eqref{eq:CR_regularised} applied to the  Bingham constitutive relation expressed via \eqref{eq:Bingham_implicit_2}. All the examples were implemented using Firedrake \cite{Rathgeber2016}. The semismooth Newton iterations are carried out until the Euclidean norm of the residual falls below $1\times 10^{-9}$, unless specified otherwise; The linear systems were solved with the sparse direct solver from MUMPS \cite{MUMPS:1}. 
Firedrake makes use of the Unified Form Language \cite{UFL}, which automatically chooses an element of the generalised derivative of the positive part $f_+$ of a function $f$ as:
\begin{equation}
\nabla f_+ = \left\{
\begin{array}{cc}
	\nabla f & \text{ if }f>0,\\
	0 & \text{ if }f\leq 0.
\end{array}
	\right.
\end{equation}

\subsection{Flow between two plates}\label{sec:flow_plates}
Consider the problem posed on the domain $\Omega = (0,4)\times (-1,1)$ with $\bm{f}=\bm{0}$. The following function solves the Bingham problem exactly (cf. \cite{Grinevich2009,Farrell2020a}):
\begin{gather}\label{eq:Bingham_Poiseuille}
\bu_e(\bm{x}):= \left\{
\begin{array}{cc}
\renewcommand{\arraystretch}{1.5}
\left(\frac{\sqrt{2}\tau_*}{4},0\right)^\top, & \textrm{ if } |x_2|\leq \frac{1}{2},\\
\left(\sqrt{2}\tau_*\left(|x_2|  - x_2^2\right),0\right)^{\top}, & \textrm{ if } |x_2|\geq \frac{1}{2},\\
\end{array}
\right.\\
p_e(\bm{x}) := \sqrt{2}\tau_*(16 - x_1).
\end{gather}
The boundary data are then chosen so as to match the given expression above.

\begin{figure}
\centering
\subfloat[C][{\centering {\normalsize Pressure obtained with the Taylor--Hood $\mathbb{P}_2$--$\mathbb{P}_1$ element.}}]{{%
	\includegraphics[width=0.95\textwidth]{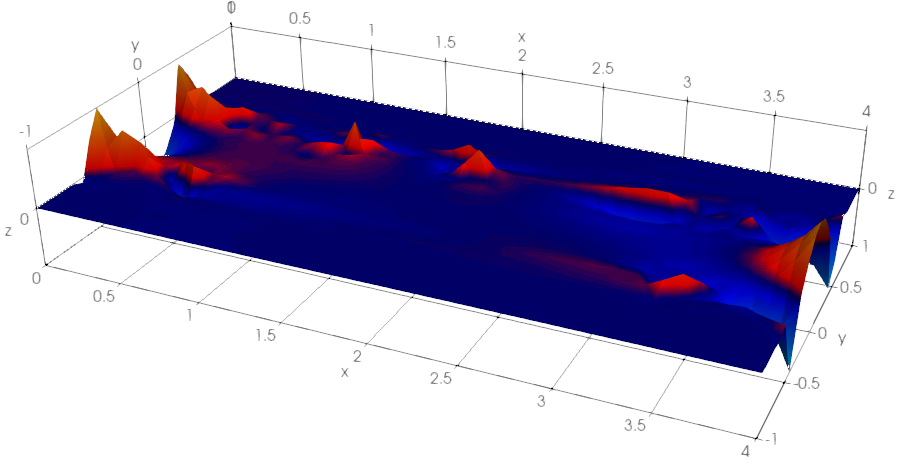}%
	}}\\
	\subfloat[D][{\centering {\normalsize Pressure obtained with the stabilised $\mathbb{P}_1$--$\mathbb{P}_1$ element.}}]{{%
	\includegraphics[width=0.95\textwidth]{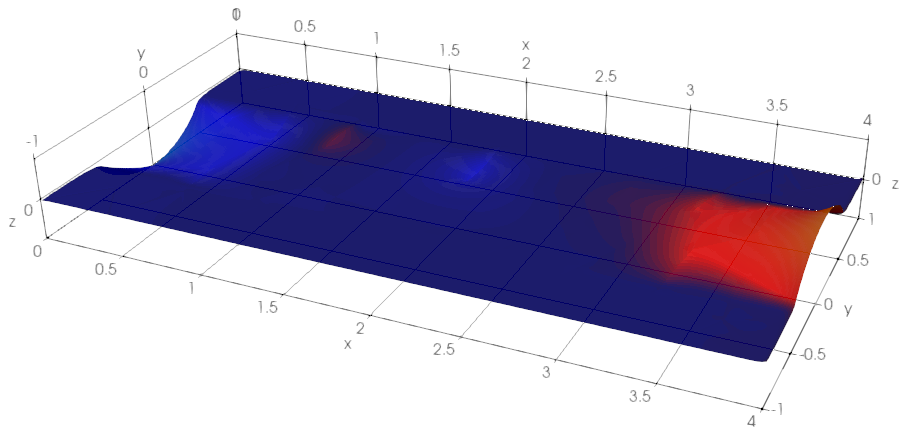}%
	}}%
	\caption{Pressure for the flow between two plates with $\varepsilon=0.001$ and $3.6\times 10^4$ dofs.}%
	\label{fig:TH_oscillation}
\end{figure}

Even though there are many families of finite element methods satisfying Assumptions \ref{as:approximability}--\ref{as:Projector_Pressure} (recall the discussion in Section \ref{sec:FEM}), which in theory lead to a convergent scheme for a big class of implicitly constituted fluids, in the particular case of viscoplastic flow it has been observed in practice that schemes for which the shear stress $\BS$ is approximated by piecewise constant functions are more appropriate (see e.g.\ \cite{Treskatis2018}). For instance, when using the Taylor--Hood $\mathbb{P}_2$--$\mathbb{P}_1$ element, oscillations that spoil the overall convergence of the nonlinear iterations might appear; Figure \ref{fig:TH_oscillation} shows e.g.\ such oscillations in the solution obtained using a mesh with $3.6\times 10^4$ degrees of freedom. For this reason, in the rest of this work we will employ instead a stabilised $\mathbb{P}_1$--$\mathbb{P}_1$ element for the velocity and pressure; i.e.\ the finite element spaces will be chosen as (recall \eqref{eq:space_discrete_stresses}):
\begin{align*}
	\Sigma^n &:= \left\{\Bsigma \in \Lsym{\infty} \colon \Bsigma|_K \in \mathbb{P}_{0}(K)^{d\times d},\, K\in\mathcal{T}_n \right\}, \\
	V^n &:= \left\{\bv \in C(\overline{\Omega})^{d} \colon \bv|_{\partial\Omega}=\bu_e,\, \bv|_K \in \mathbb{P}_{1}(K)^d,\, K\in\mathcal{T}_n \right\}, \\
M^n &:= \left\{q \in \Lmean{\infty} \colon q|_K \in \mathbb{P}_{1}(K),\, K\in\mathcal{T}_n \right\},
\end{align*}
and the following stabilisation term is added to the right-hand-side of the discrete equation \eqref{eq:discrete_mass}:
\begin{equation}
 0.2 h^2 \int_\Omega \nabla p^n\cdot \nabla q.
\end{equation}
\txtbl{In all the numerical experiments, the phrase ``degrees of freedom'' (or ``\#dofs'') refers to the total number of unknowns associated to the stress, velocity, and pressure.}

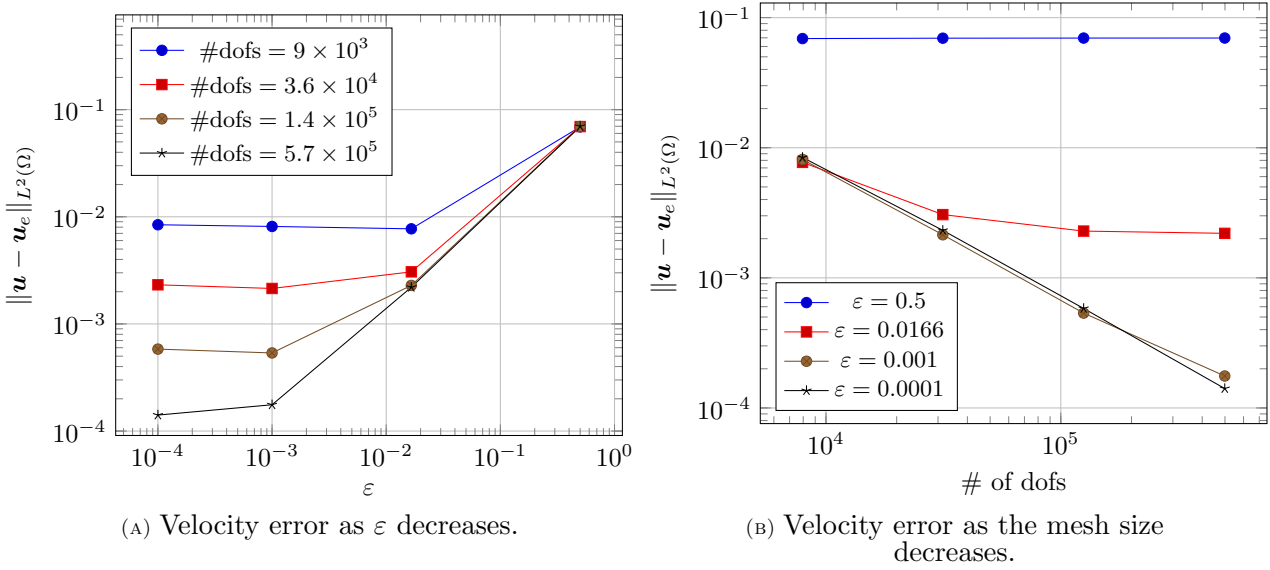
\begin{figure}
\centering
\resizebox{\linewidth}{!}{
	\subfloat[C][{\centering {\large Velocity error as $\varepsilon$ decreases.}}]{{%
\begin{tikzpicture}
	\begin{loglogaxis}[
		xlabel={$\varepsilon$},
		ylabel={$\|\bm{u} - \bm{u}_e\|_{L^2(\Omega)}$},
		enlarge y limits={abs=11},
		ymin=0.001,
		normalsize,
		legend style={font=\small},
		legend pos=north west,
		grid=major,
		legend entries={$\mathrm{\# dofs}=9\times 10^3$,$\mathrm{\# dofs}=3.6\times 10^4$,$\mathrm{\# dofs}=1.4\times 10^5$,$\mathrm{\# dofs}=5.7\times 10^5$},
	]
		\addplot table {data_errors_new/vel_epserror_nref_1_MAX.dat};
		\addplot table {data_errors_new/vel_epserror_nref_2_MAX.dat};
		\addplot table {data_errors_new/vel_epserror_nref_3_MAX.dat};
		\addplot table {data_errors_new/vel_epserror_nref_4_MAX.dat};
	\end{loglogaxis}
\end{tikzpicture}%
}}%
\subfloat[D][{\centering {\large Velocity error as the mesh size decreases.}}]{{%
\begin{tikzpicture}
	\begin{loglogaxis}[
		xlabel={\# of dofs},
		ylabel={$\|\bm{u} - \bm{u}_e\|_{L^2(\Omega)}$},
		legend style={font=\small},
		legend pos=south west,
		normalsize,
		grid=major,
		legend entries={$\varepsilon=0.5$,$\varepsilon=0.0166$,$\varepsilon=0.001$,$\varepsilon=0.0001$},
	]
		\addplot table {data_errors_new/vel_dofserror_eps_0.5_MAX.dat};
		\addplot table {data_errors_new/vel_dofserror_eps_0.0166_MAX.dat};
		\addplot table {data_errors_new/vel_dofserror_eps_0.001_MAX.dat};
		\addplot table {data_errors_new/vel_dofserror_eps_1e-4_MAX.dat};
	\end{loglogaxis}
\end{tikzpicture}%
	}}%
}%
	\caption{$L^2$ error of the velocity for the Bingham flow between two plates.}%
	\label{fig:vel_error}
\end{figure}

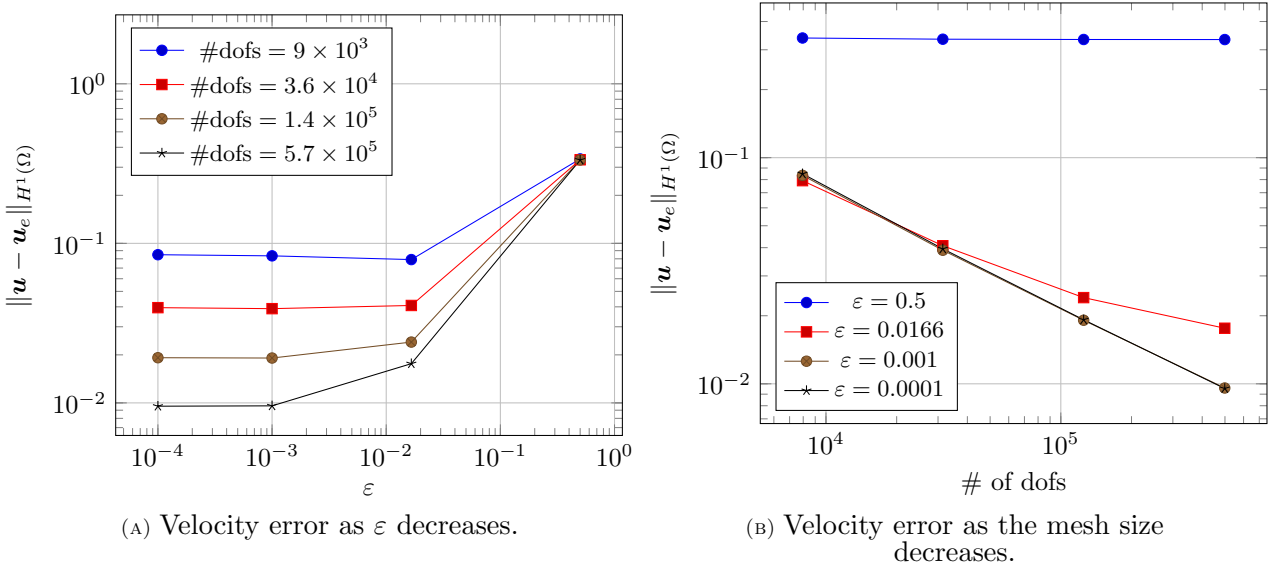
\begin{figure}
\centering
\resizebox{\linewidth}{!}{
	\subfloat[C][{\centering {\large Velocity error as $\varepsilon$ decreases.}}]{{%
\begin{tikzpicture}
	\begin{loglogaxis}[
		xlabel={$\varepsilon$},
		ylabel={$\|\bm{u} - \bm{u}_e\|_{H^1(\Omega)}$},
		enlarge y limits={abs=8},
		ymin=0.05,
		normalsize,
		legend style={font=\small},
		legend pos=north west,
		grid=major,
		legend entries={$\mathrm{\# dofs}=9\times 10^3$,$\mathrm{\# dofs}=3.6\times 10^4$,$\mathrm{\# dofs}=1.4\times 10^5$,$\mathrm{\# dofs}=5.7\times 10^5$},
	]
		\addplot table {data_errors_new/velgrad_epserror_nref_1_MAX.dat};
		\addplot table {data_errors_new/velgrad_epserror_nref_2_MAX.dat};
		\addplot table {data_errors_new/velgrad_epserror_nref_3_MAX.dat};
		\addplot table {data_errors_new/velgrad_epserror_nref_4_MAX.dat};
	\end{loglogaxis}
\end{tikzpicture}%
}}%
\subfloat[D][{\centering {\large Velocity error as the mesh size decreases.}}]{{%
\begin{tikzpicture}
	\begin{loglogaxis}[
		xlabel={\# of dofs},
		ylabel={$\|\bm{u} - \bm{u}_e\|_{H^1(\Omega)}$},
		legend style={font=\small},
		legend pos=south west,
		normalsize,
		grid=major,
		legend entries={$\varepsilon=0.5$,$\varepsilon=0.0166$,$\varepsilon=0.001$,$\varepsilon=0.0001$},
	]
		\addplot table {data_errors_new/velgrad_dofserror_eps_0.5_MAX.dat};
		\addplot table {data_errors_new/velgrad_dofserror_eps_0.0166_MAX.dat};
		\addplot table {data_errors_new/velgrad_dofserror_eps_0.001_MAX.dat};
		\addplot table {data_errors_new/velgrad_dofserror_eps_1e-4_MAX.dat};
	\end{loglogaxis}
\end{tikzpicture}%
	}}%
}%
	\caption{$H^1$ error of the velocity for the Bingham flow between two plates.}%
	\label{fig:velgrad_error}
\end{figure}

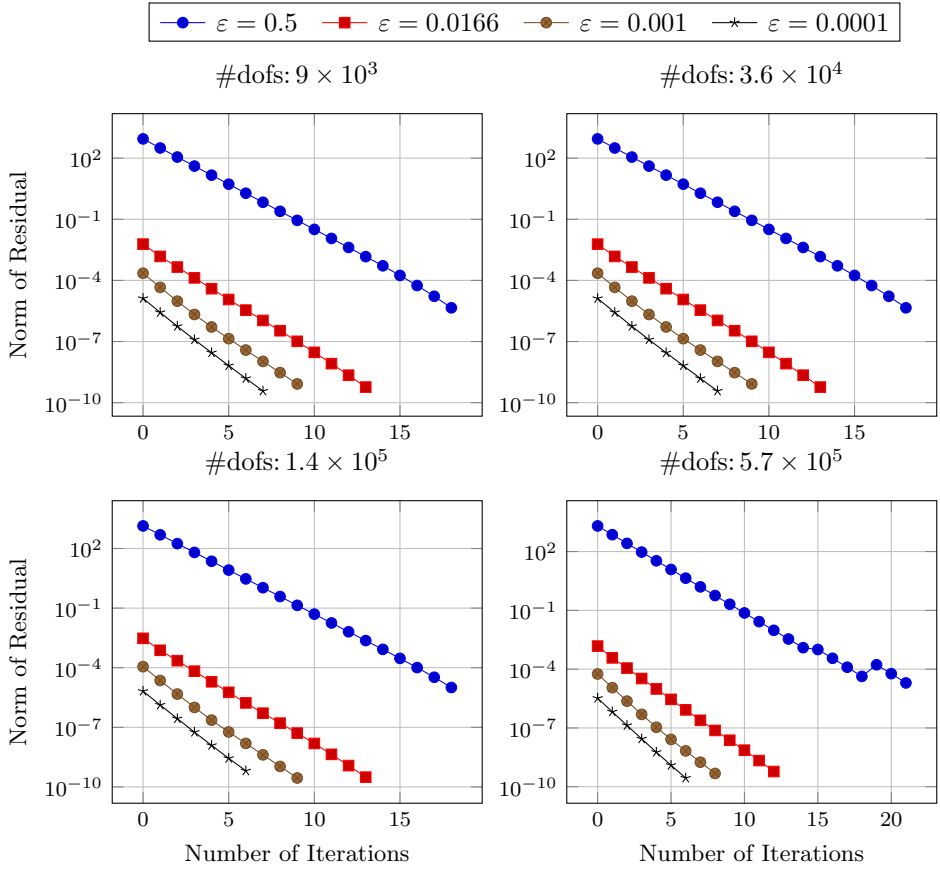
\begin{figure}
	\centering
\begin{tikzpicture}
    \begin{groupplot}[group style={
                      group name=myplot,
		      vertical sep=32pt,
		      horizontal sep=32pt,
                      group size= 2 by 2},height=5.5cm,width=6.5cm]
		      \nextgroupplot[small, title=\#dofs${:}\, 9\times 10^3$,ymode=log,grid=major,ylabel={Norm of Residual }]
		\addplot table {data_residual/residual_nref_2_eps_0.5.dat};\label{plots:eps1}
		\addplot table {data_residual/residual_nref_2_eps_0.0166.dat};\label{plots:eps2}
		\addplot table {data_residual/residual_nref_2_eps_0.001.dat};\label{plots:eps3}
		\addplot table {data_residual/residual_nref_2_eps_0.0001.dat};\label{plots:eps4}

        \nextgroupplot[title=\#dofs${:}\, 3.6\times 10^4$,small,ymode=log,grid=major]
		\addplot table {data_residual/residual_nref_2_eps_0.5.dat};
		\addplot table {data_residual/residual_nref_2_eps_0.0166.dat};
		\addplot table {data_residual/residual_nref_2_eps_0.001.dat};
		\addplot table {data_residual/residual_nref_2_eps_0.0001.dat};
        \nextgroupplot[title=\#dofs${:}\, 1.4\times 10^5$,small,ymode=log,grid=major,xlabel={Number of Iterations},ylabel={Norm of Residual }]
		\addplot table {data_residual/residual_nref_3_eps_0.5.dat};
		\addplot table {data_residual/residual_nref_3_eps_0.0166.dat};
		\addplot table {data_residual/residual_nref_3_eps_0.001.dat};
		\addplot table {data_residual/residual_nref_3_eps_0.0001.dat};
        \nextgroupplot[title=\#dofs${:}\, 5.7\times 10^5$,small,ymode=log,grid=major,xlabel={Number of Iterations}]
		\addplot table {data_residual/residual_nref_4_eps_0.5.dat};
		\addplot table {data_residual/residual_nref_4_eps_0.0166.dat};
		\addplot table {data_residual/residual_nref_4_eps_0.001.dat};
		\addplot table {data_residual/residual_nref_4_eps_0.0001.dat};
    \end{groupplot}
\path (myplot c1r1.north west|-current bounding box.north)--
      coordinate(legendpos)
      (myplot c2r1.north east|-current bounding box.north);
\matrix[
    matrix of nodes,
    anchor=south,
    draw,
    inner sep=0.2em,
    draw
  ]at([yshift=1ex]legendpos)
  {
    \ref{plots:eps1}& $\varepsilon = 0.5$&[5pt]
    \ref{plots:eps2}& $\varepsilon = 0.0166$ &[5pt]
    \ref{plots:eps3}& $\varepsilon = 0.001$ &[5pt]
    \ref{plots:eps4}& $\varepsilon = 0.0001$\\};
\end{tikzpicture}
\caption{Residual norm vs. number of (semismooth) Newton iterations for the Bingham flow between two plates, for different mesh sizes.}%
\label{fig:iterations_dofs}
\end{figure}

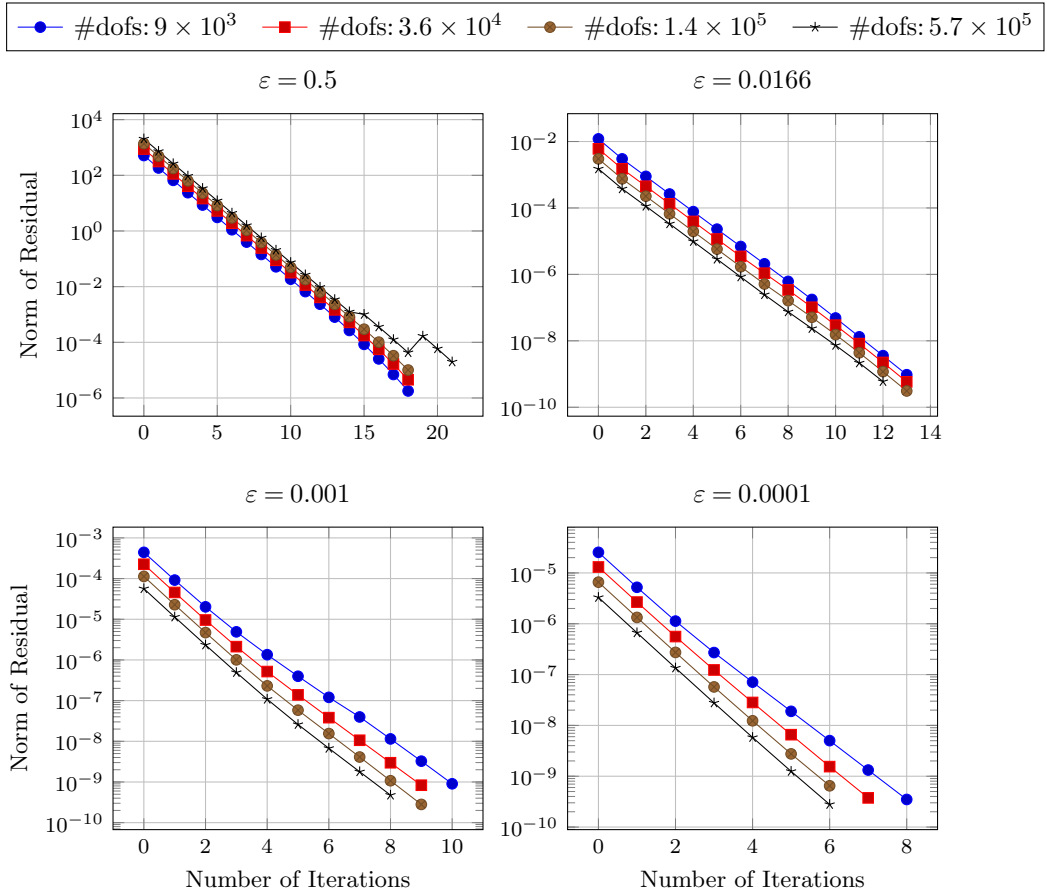
\begin{figure}
\begin{tikzpicture}
    \begin{groupplot}[group style={
                      group name=myplot,
		      vertical sep=42pt,
		      horizontal sep=32pt,
                      group size= 2 by 2},height=5.5cm,width=6.5cm]
		      \nextgroupplot[small, title=$\varepsilon\:{=}\: 0.5$,ymode=log,grid=major,ylabel={Norm of Residual }]
		\addplot table {data_residual/residual_nref_1_eps_0.5.dat};\label{plots:nref1}
		\addplot table {data_residual/residual_nref_2_eps_0.5.dat};\label{plots:nref2}
		\addplot table {data_residual/residual_nref_3_eps_0.5.dat};\label{plots:nref3}
		\addplot table {data_residual/residual_nref_4_eps_0.5.dat};\label{plots:nref4}

			\nextgroupplot[title=$\varepsilon\:{=}\: 0.0166$,small,ymode=log,grid=major]
		\addplot table {data_residual/residual_nref_1_eps_0.0166.dat};
		\addplot table {data_residual/residual_nref_2_eps_0.0166.dat};
		\addplot table {data_residual/residual_nref_3_eps_0.0166.dat};
		\addplot table {data_residual/residual_nref_4_eps_0.0166.dat};
		\nextgroupplot[title=$\varepsilon\:{=}\: 0.001$,small,ymode=log,grid=major,xlabel={Number of Iterations},ylabel={Norm of Residual }]
		\addplot table {data_residual/residual_nref_1_eps_0.001.dat};
		\addplot table {data_residual/residual_nref_2_eps_0.001.dat};
		\addplot table {data_residual/residual_nref_3_eps_0.001.dat};
		\addplot table {data_residual/residual_nref_4_eps_0.001.dat};
		\nextgroupplot[title=$\varepsilon\:{=}\: 0.0001$,small,ymode=log,grid=major,xlabel={Number of Iterations}]
		\addplot table {data_residual/residual_nref_1_eps_0.0001.dat};
		\addplot table {data_residual/residual_nref_2_eps_0.0001.dat};
		\addplot table {data_residual/residual_nref_3_eps_0.0001.dat};
		\addplot table {data_residual/residual_nref_4_eps_0.0001.dat};
    \end{groupplot}
\path (myplot c1r1.north west|-current bounding box.north)--
      coordinate(legendpos)
      (myplot c2r1.north east|-current bounding box.north);
\matrix[
    matrix of nodes,
    anchor=south,
    draw,
    inner sep=0.2em,
    draw
  ]at([yshift=1ex]legendpos)
  {
	  \ref{plots:nref1}& \#dofs${:}\, 9\times 10^3$ &[5pt]
    \ref{plots:nref2}& \#dofs${:}\, 3.6\times 10^4$ &[5pt]
    \ref{plots:nref3}& \#dofs${:}\, 1.4\times 10^5$ &[5pt]
    \ref{plots:nref4}& \#dofs${:}\, 5.7\times 10^5$\\};
\end{tikzpicture}
\caption{Residual norm vs.\ number of (semismooth) Newton iterations for the Bingham flow between two plates, for different values of $\varepsilon$.}%
\label{fig:iterations_eps}
\end{figure}

In this example we set $\tau_*=1$, and continuation was employed to obtain better initial guesses for the nonlinear iterations; in particular, the problem was solved for the values $\varepsilon\in \{0.5,0.0166,0.001,0.0001\}$, with each solution being used as the initial guess in the iterations, for the next smaller value of $\varepsilon$. \txtbl{Figures \ref{fig:vel_error} and \ref{fig:velgrad_error} show the $L^2$ and $H^1$ error of the velocity as $\varepsilon$ and the mesh size decrease, respectively}; it can be observed that for values less than $\varepsilon=0.001$, it is the discretisation error that dominates. Figures \ref{fig:iterations_dofs} and \ref{fig:iterations_eps} show the decrease in the residual norm as the number of semismooth Newton iterations increases. The different values of $\varepsilon$ and the mesh size result in very similar behaviour for the residual, and in particular we see that the solver does not seem to exhibit any mesh-dependence.

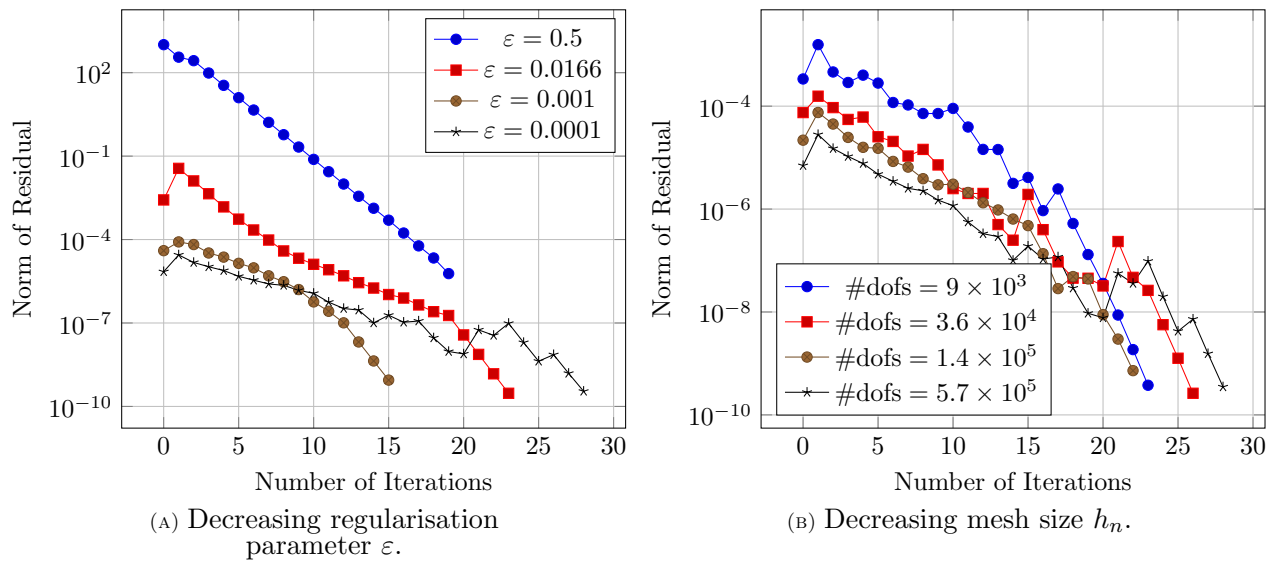
\begin{figure}
\centering
\resizebox{\linewidth}{!}{
	\subfloat[C][{\centering {\large Decreasing regularisation parameter $\varepsilon$.}}]{{%
\begin{tikzpicture}
	\begin{semilogyaxis}[
		xlabel={Number of Iterations},
		ylabel={Norm of Residual},
		grid=major,
		legend entries={$\varepsilon=0.5$,$\varepsilon=0.0166$,$\varepsilon=0.001$,$\varepsilon=0.0001$},
	]
		\addplot table {data_residual/residual_nref_4_eps_0.5_MAX.dat};
		\addplot table {data_residual/residual_nref_4_eps_0.0166_MAX.dat};
		\addplot table {data_residual/residual_nref_4_eps_0.001_MAX.dat};
		\addplot table {data_residual/residual_nref_4_eps_0.0001_MAX.dat};
	\end{semilogyaxis}
\end{tikzpicture}
}}%
\subfloat[D][{\centering {\large Decreasing mesh size $h_n$.}}]{{%
\begin{tikzpicture}
	\begin{semilogyaxis}[
		xlabel={Number of Iterations},
		ylabel={Norm of Residual},
		legend pos=south west,
		grid=major,
		legend entries={$\mathrm{\# dofs}=9\times 10^3$,$\mathrm{\# dofs}=3.6\times 10^4$,$\mathrm{\# dofs}=1.4\times 10^5$,$\mathrm{\# dofs}=5.7\times 10^5$},
	]
		\addplot table {data_residual/residual_nref_1_eps_0.0001_MAX.dat};
		\addplot table {data_residual/residual_nref_2_eps_0.0001_MAX.dat};
		\addplot table {data_residual/residual_nref_3_eps_0.0001_MAX.dat};
		\addplot table {data_residual/residual_nref_4_eps_0.0001_MAX.dat};
	\end{semilogyaxis}
\end{tikzpicture}%
	}}%
}%
\caption{Residual norm vs.\ number of (semismooth) Newton iterations for the Bingham flow between two plates using \eqref{eq:Bingham_implicit_1}.}%
	\label{fig:iterations_MAX}
\end{figure}

For the sake of comparison, Figure \ref{fig:iterations_MAX} shows the decrease in the residual norm when using the alternative expression \eqref{eq:Bingham_implicit_1} for the constitutive relation; recall that, as seen in Section \ref{sec:Bingham}, this expression does not define a semismooth function. To ensure convergence of the nonlinear iterations, in this case we had to employ a continuation scheme with respect to $\varepsilon$: given two previously computed solutions $\bm{z}_1,\bm{z}_2$, corresponding to the values $\varepsilon_1,\varepsilon_2$, respectively, the initial guess for the semismooth Newton iteration for $\varepsilon$ is chosen as
\begin{equation}
	\frac{\varepsilon - \varepsilon_2}{\varepsilon_2-\varepsilon_1}(\bm{z}_2 - \bm{z}_1) + \bm{z}_2.
\end{equation}
Although convergence is achieved, the behaviour is somewhat more erratic in comparison to that obtained using \eqref{eq:Bingham_implicit_2}, and it takes a higher number of iterations and continuation steps for the solver to converge. This supports the results from Section \ref{sec:Bingham}, which suggest that the expression \eqref{eq:Bingham_implicit_2} should be more advantageous.

\subsection{Lid-driven cavity}

The lid-driven cavity problem has been used as a test for many numerical schemes throughout the years (see e.g.\ \cite{Syrakos2013,Aposporidis2011}). In this example we solve the unsteady problem on the  domain defined by the square $\Omega = (0,1)^2$; the following boundary conditions are imposed on the velocity field for all times $t\in[0,T]$:

\begin{align}
\partial\Omega_1 &:= (0,1)\times \{1\}, \qquad &\partial\Omega_2 := \partial\Omega \setminus \partial\Omega_1,\notag \\
\bm{u} &= \bm{0}\qquad &\text{ on }\partial\Omega_2, \label{eq:lid_BCs}\\
\bm{u} &= (1,0)^\text{T}\qquad &\text{ on }\partial\Omega_1.\notag
\end{align}

The time derivative is discretised using implicit Euler's method with a time step $\Delta t$, and the goal is to compute the steady state of the system. In this example we use an absolute tolerance for the nonlinear iterations of $10^{-7}$, and we deem the solution to have reached a steady state when the difference between two consecutive solutions $\|\bu^{k+1} - \bu^k\|_{L^2(\Omega)}$ is less than $10^{-6}$.

\begin{figure}
\centering
\subfloat[C][{\centering {\normalsize $\tau_*=3$; $\Delta t = 0.0005$}.}]{{%
	\includegraphics[width=0.5\textwidth]{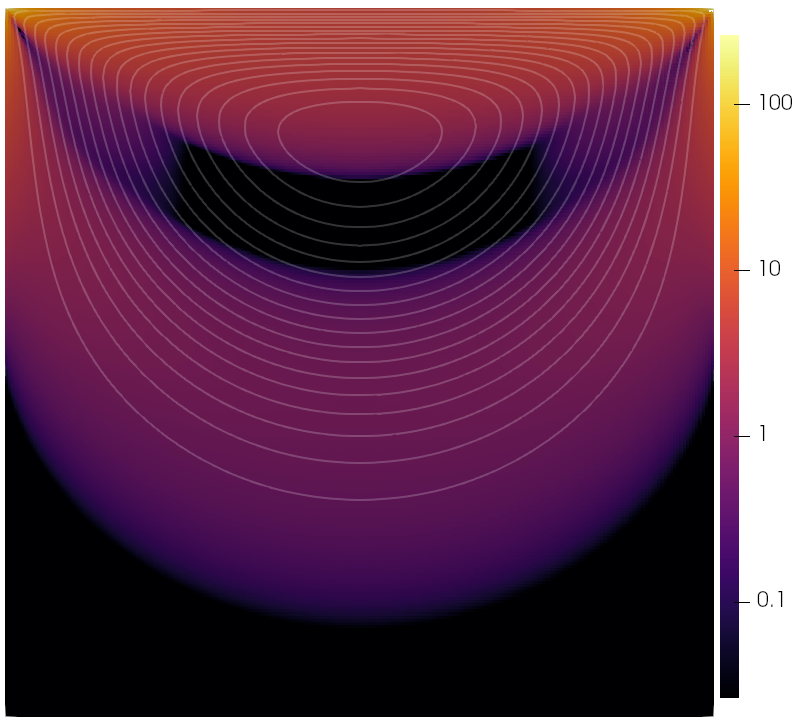}%
	}}%
	\subfloat[D][{\centering {\normalsize $\tau_*=3$; $\Delta t = 0.0005$}.}]{{%
	\includegraphics[width=0.5\textwidth]{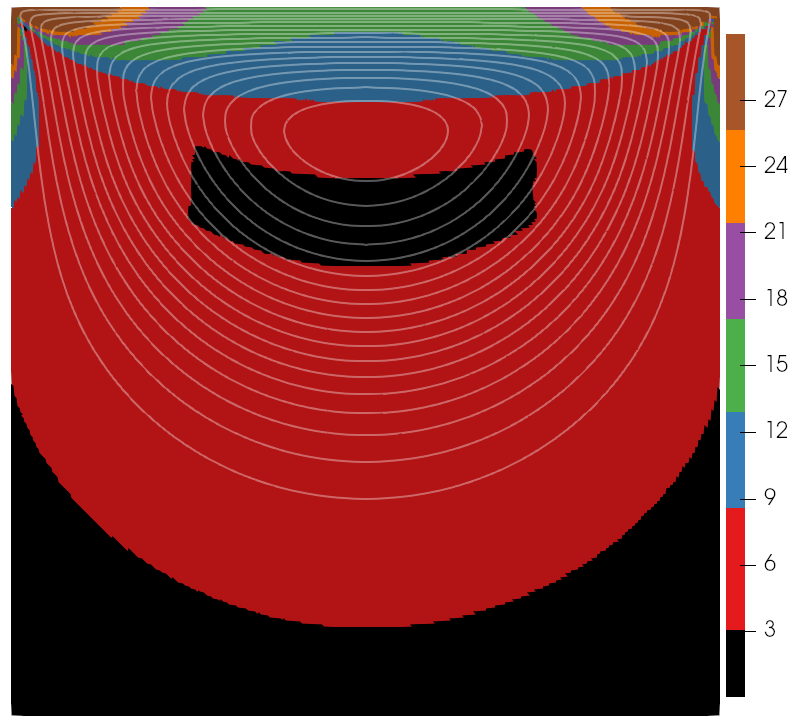}%
	}}\\
	\subfloat[A][{\centering {\normalsize $\tau_*=20$; $\Delta t = 0.0001$}.}]{{%
	\includegraphics[width=0.5\textwidth]{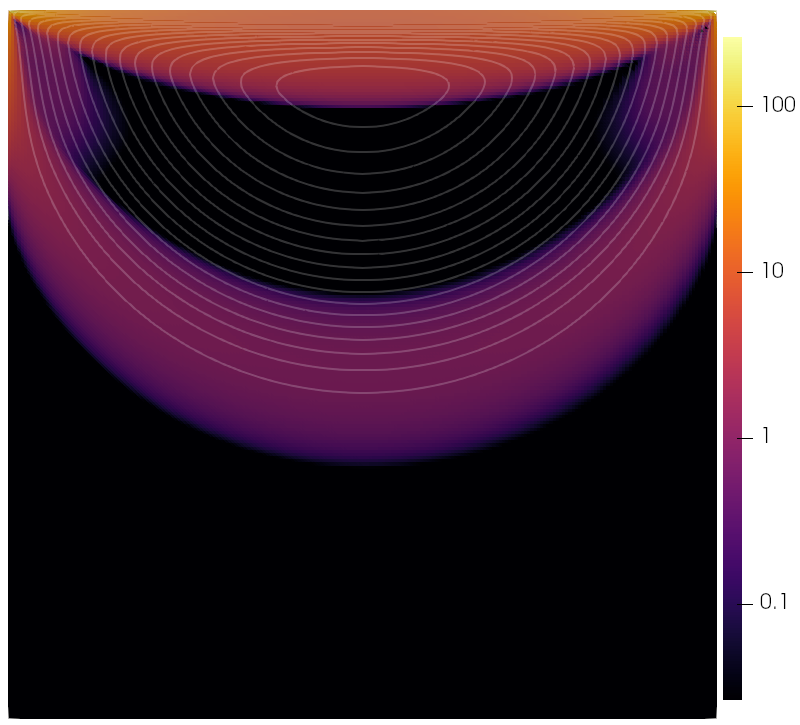}%
	}}%
	\subfloat[B][{\centering {\normalsize $\tau_*=20$; $\Delta t = 0.0001$}.}]{{%
	\includegraphics[width=0.5\textwidth]{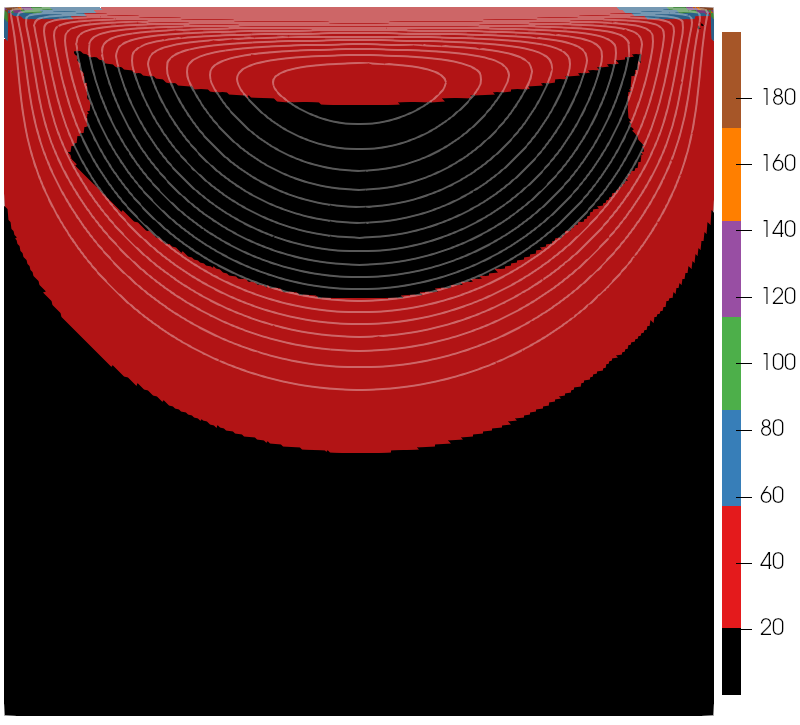}%
	}}\\
	\caption{Magnitude of $\Du$ (left) and $\BS$ (right) for the Bingham cavity with $\varepsilon=10^{-4}$ and $\tau_*\in \{3,20\}$.}%
	\label{fig:Cavity}
\end{figure}

\txtbl{Figure \ref{fig:Cavity} shows comparative plots of the magnitude of the symmetric velocity gradient $\Du$ and the stress $\BS$ of the steady state, for two different values of the yield stress, obtained using a discretisation with $4.5\times 10^5$ degrees of freedom. The scheme yields a solution that exhibits the expected behaviour, i.e.\ a plug region appears in the upper half of the cavity, and it increases in size as $\tau_*$ increases. Table \ref{tb:cavity_} shows, for several values of the yield stress $\tau_*$, the maximum of the stream function $\psi_{\mathrm{max}}$, and the vertical coordinate of the center of the vortex that appears in the middle of the cavity. The obtained values are in good agreement with \cite{Syrakos2013}.}


\begin{table}
\centering
\captionsetup{justification=centering}
\begin{tabular}{c c c c c c c}
\toprule
$\tau_*$ & 0.5 & 3 & 5 & 10 & 20 & 40 \\
\midrule
$\psi_{\mathrm{max}}$ & 0.09472 & 0.07716 & 0.06916 & 0.05734  & 0.04554 & 0.03473 \\
$y_c$ & 0.7768 & 0.8170 & 0.8348 & 0.8616 & 0.8884 & 0.9107  \\
  \bottomrule
\end{tabular}\\
\caption{\label{tb:cavity_}
Strength and vertical coordinate of the central vortex for the lid-driven cavity problem}
\end{table}


\subsection{Expansion-contraction channel}
{\color{black}
	In this section we consider the flow of a Bingham fluid on a channel whose width expands and then contracts; the geometry of the problem is depicted in Figure \ref{fig:exp_cont_channel}. The problem has already been non-dimensionalised, and in the current setting the Bingham number $\mathrm{Bn}$ takes the role of the yield stress; for details about the non-dimensionalisation procedure see \cite{Marly2017}, where this problem was solved employing the augmented Lagrangian method and a finite difference discretisation.

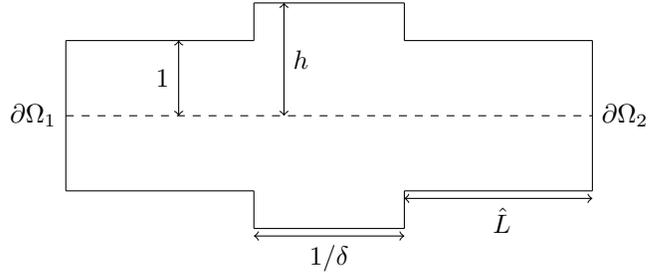
\begin{figure}
\begin{tikzpicture}
	\draw (-3.5,-1) -- (-3.5,1);
	\draw (-3.5,1) -- (-1,1);
	\draw (-1,1) -- (-1,1.5);
	\draw (-1,1.5) -- (1,1.5);
	\draw (1,1.5) -- (1,1);
	\draw (1,1) -- (3.5,1);
	\draw (3.5,1) -- (3.5,-1);
	\draw (3.5,-1) -- (1,-1);
	\draw (1,-1) -- (1,-1.5);
	\draw (1,-1.5) -- (-1,-1.5);
	\draw (-1,-1.5) -- (-1,-1);
	\draw (-1,-1) -- (-3.5,-1);
	\node [left] at (-3.5,0) {$\partial\Omega_1$};
	\node [right] at (3.5,0) {$\partial\Omega_2$};
	\draw [<->] (-2,0) -- (-2,1);
	\node [left] at (-2,0.5) {$1$};
	\draw [<->] (-0.6,0) -- (-0.6,1.5);
	\node [right] at (-0.6,0.75) {$h$};
	\draw [dashed] (-3.5,0) -- (3.5,0);
	\draw [<->] (-1,-1.6) -- (1,-1.6);
	\node [below] at (0, -1.6) {$1/\delta$};
	\draw [<->] (1,-1.1) -- (3.5,-1.1);
	\node [below] at (2.3, -1.1) {$\hat{L}$};
\end{tikzpicture}%
	\caption{Geometry of the expansion-contraction problem.}%
	\label{fig:exp_cont_channel}
\end{figure}

The boundary condition for the velocity at the inflow and outflow boundaries $\partial\Omega_1\cup \partial\Omega_2$ is chosen to be a fully developed Poiseuille flow (cf. Section \ref{sec:flow_plates}). On the rest of the boundary $\partial\Omega \setminus (\partial\Omega_1 \cup\partial\Omega_2)$ we impose a no-slip condition $\bu=0$. Figure \ref{fig:exp_cont_short} shows the magnitude of $\Du$ (using a logarithmic scale for the color map) for the solution obtained using a mesh with 330819 degrees of freedom, and the values $\hat{L}=4$, $\delta=\frac{1}{2}$, $h=2$. Figure \ref{fig:exp_cont_pressure} then shows the variation of the pressure in the horizontal direction; the pressure exhibits a linear variation before and after the cavity, as expected from a Poiseuille profile.

\begin{figure}
\centering
\subfloat[C][{\centering {\normalsize $\mathrm{Bn}=2$}}]{{%
	\includegraphics[width=0.80\textwidth,height=32ex]{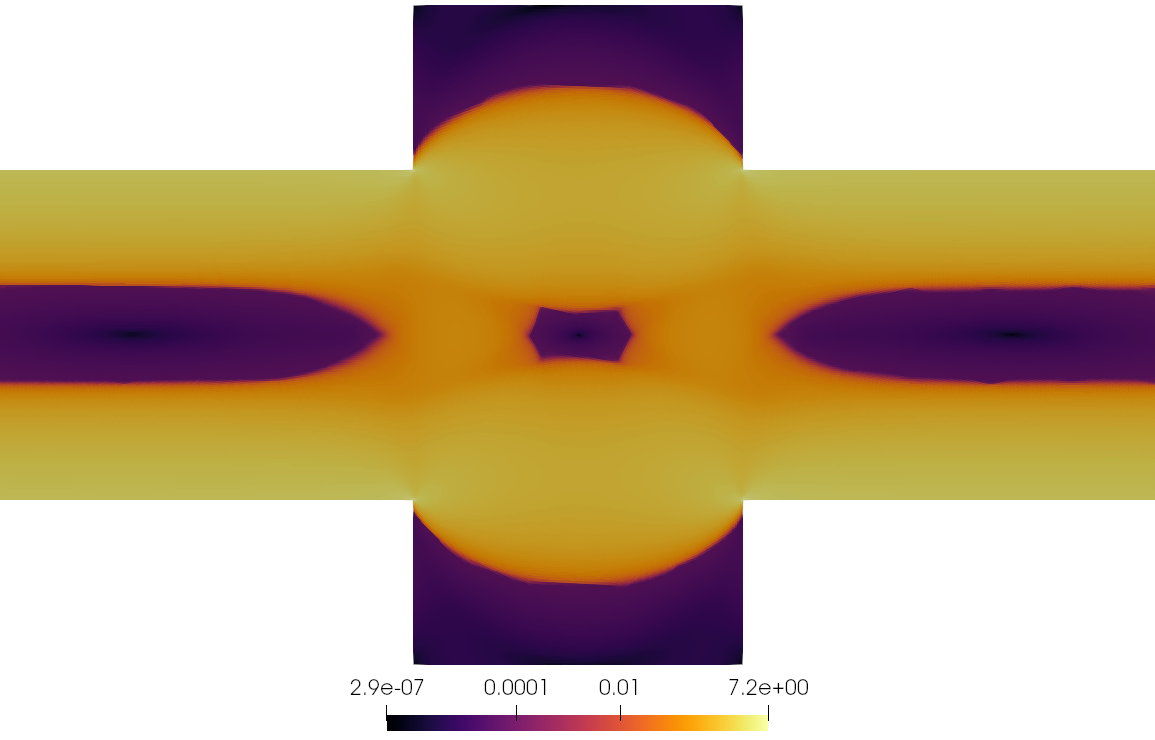}%
	}}\\
	\subfloat[D][{\centering {\normalsize $\mathrm{Bn}=10$}}]{{%
	\includegraphics[width=0.80\textwidth,height=32ex]{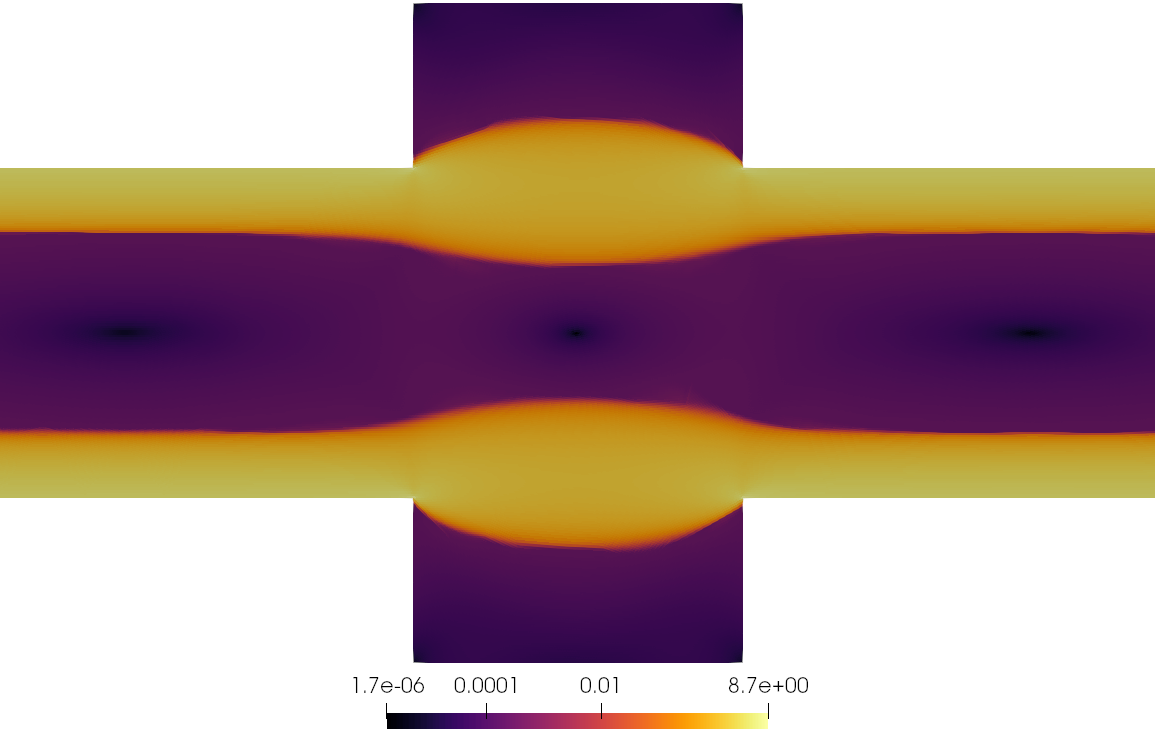}%
	}}%
	\caption{Magnitude of $|\Du|$ for the expansion-contraction problem with $\varepsilon=10^{-4}$, $\hat{L}=4$, $\delta=\frac{1}{2}$, $h=2$.}%
	\label{fig:exp_cont_short}
\end{figure}

\begin{figure}
	\centering
	\includegraphics[width=0.70\textwidth]{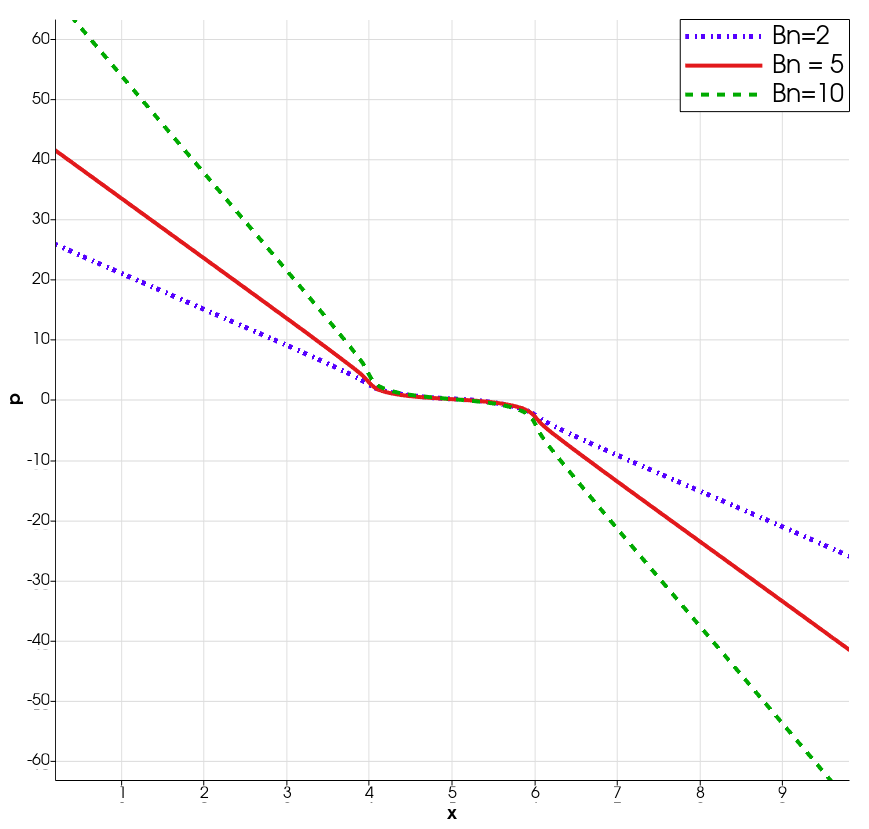}%
	\caption{Pressure variation along a horizontal line with a height halfway between the boundary of the channel and the beginning of the Poiseuille plug.}%
	\label{fig:exp_cont_pressure}
\end{figure}

As the height $h$ gets closer to $1$, the plug in the cavity breaks into two dead zones; Figure \ref{fig:exp_cont_long} shows two different solutions that exhibit this feature. Table \ref{tb:exp_cont} shows the (non-dimensionalised) length of the dead zone for several values of the Bingham number. Although the interfaces between yielded and unyielded regions is not as sharp as those obtained with the augmented Lagrangian method in \cite{Marly2017}, we find that the proposed algorithm is still capable of producing a behaviour consistent with \cite{Marly2017}, including e.g. the small plug regions that appear as the channel expands or contracts (cf.\ Figure \ref{fig:exp_cont_long} \textsc{(b)}).

\begin{table}
\centering
\captionsetup{justification=centering}
\begin{tabular}{c c c c c c c c }
\toprule
$\mathrm{Bn}$ & 2 & 5 & 10 & 20 & 30 & 40 & 50 \\
\midrule
$L_d$ & 0.0469 & 0.06453 & 0.07813 & 0.09688  & 0.10626 & 0.11876 & 0.12500 \\
  \bottomrule
\end{tabular}\\
\caption{\label{tb:exp_cont}
Length of the dead zone for the expansion-contraction problem with $\hat{L}=3$, $\delta=\frac{1}{5}$, $h=\frac{6}{5}$, for several values of the Bingham number.}
\end{table}

\begin{figure}
\centering
\subfloat[C][{\centering {\normalsize $\mathrm{Bn}=5$}}]{{%
	\includegraphics[width=0.80\textwidth,height=32ex]{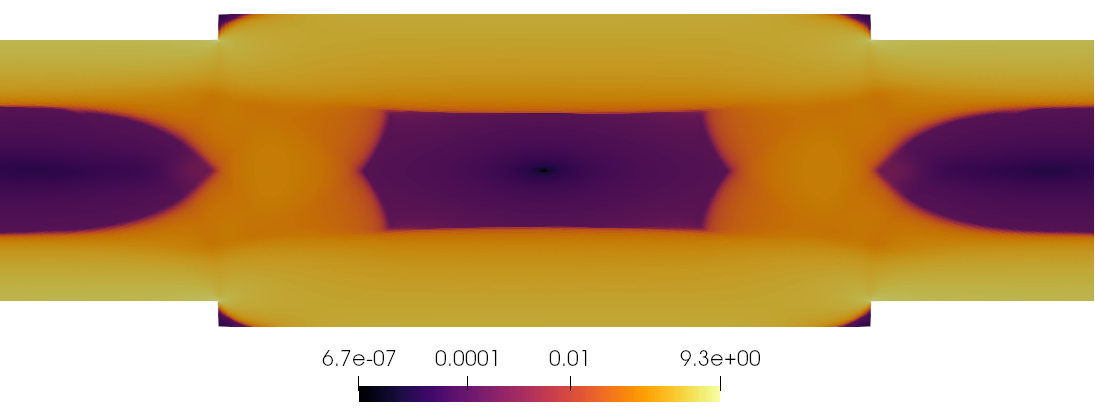}%
	}}\\
	\subfloat[D][{\centering {\normalsize $\mathrm{Bn}=50$}}]{{%
	\includegraphics[width=0.80\textwidth,height=32ex]{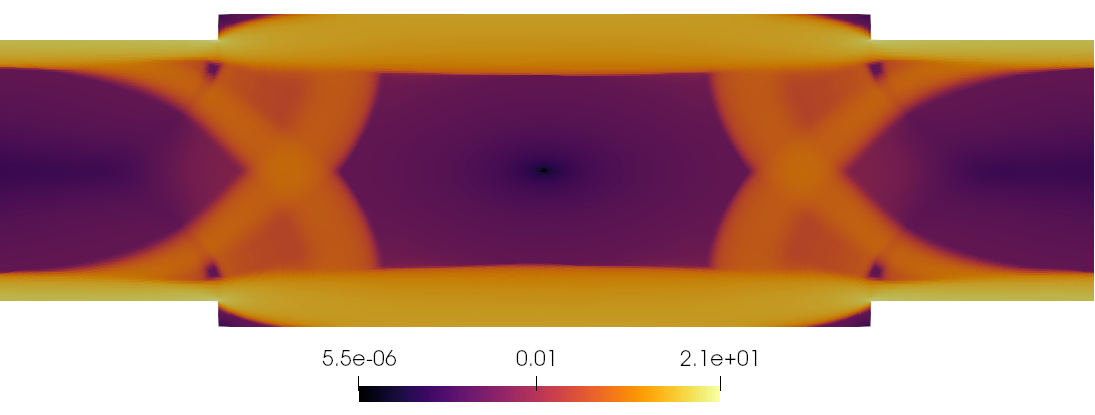}%
	}}%
	\caption{Magnitude of $|\Du|$ for the expansion-contraction problem with $\varepsilon=10^{-4}$, $\hat{L}=3$, $\delta=\frac{1}{5}$, $h=\frac{6}{5}$.}%
	\label{fig:exp_cont_long}
\end{figure}

}
\subsection*{Acknowledgements}
\begin{acknowledgement}
The author would like to thank T.\ Surowiec for helpful pointers.
\end{acknowledgement}
\bibliographystyle{plain}
\bibliography{bibliography}
\end{document}